\AddToHook{begindocument/before}{} 
\DeclareHookRule{begindocument}{hyperref}{before}{my-mdpi} 

\documentclass[preprints,article,accept,moreauthors,pdftex]{my-mdpi}

\makeatletter
\renewcommand{\maketag@@@}[1]{\hbox{\m@th\normalsize\normalfont#1}}
\makeatother

\usepackage{graphicx}
\usepackage{subcaption}
\usepackage{subfloat}
\usepackage{environ}
\usepackage{amsmath,amssymb}

\firstpage{1} 
\pubvolume{15}
\issuenum{10}
\articlenumber{1911}
\pubyear{2023}
\copyrightyear{2023}
\externaleditor{\vspace{-12pt}Academic Editor: Yongshun Liang, Saurabh Verma and Arulprakash Gowrisankar}
\datereceived{25 August 2023} 
\daterevised{7 October 2023} 
\dateaccepted{10 October 2023} 
\datepublished{12 October 2023} 
\hreflink{https://doi.org/10.3390/sym15101911} 
\doinum{10.3390/sym15101911}

\pdfoutput=1 

\usepackage{dsfont} 

\Title{On Sharp Bounds of Local Fractional Metric Dimension for Certain Symmetrical Algebraic Structure Graphs}

\TitleCitation{On Sharp Bounds of Local Fractional Metric Dimension for Certain Symmetrical Algebraic Structure Graphs}

\Author{Amal S. Alali $^{1}$\orcidA{}, 
Shahbaz Ali $^{2}$\orcidB{}, 
Muhammad Adnan $^{2}$ 
and Delfim F. M. Torres $^{3,}$*\orcidD{}}

\AuthorNames{Amal S. Alali, Shahbaz Ali, Muhammad Adnan and Delfim F. M. Torres}

\AuthorCitation{Alali,~A.S.; Ali,~S.; Adnan,~M.; Torres,~D.F.M.}

\address{%
$^{1}$ \quad Department of Mathematical Sciences, College of Science, 
Princess Nourah bint Abdulrahman University, P.O.~Box 84428, Riyadh 11671, Saudi Arabia; asalali@pnu.edu.sa\\
$^{2}$ \quad Department of Mathematics, The Islamia University of Bahawalpur, 
Rahim Yar Kahn Campus, \mbox{Rahim Yar Khan} 64200, Pakistan; 
shahbazali@iub.edu.pk (S.A.); s21rmath3e01052@iub.edu.pk (M.A.)\\
$^{3}$ \quad Center for Research and Development in Mathematics and Applications (CIDMA),
\mbox{Department of Mathematics}, University of Aveiro, 3810-193 Aveiro, Portugal}

\corres{Correspondence: delfim@ua.pt}

\abstract{The smallest set of vertices needed 
to differentiate or categorize every other vertex in a graph 
is referred to as the graph's metric dimension. 
Finding the class of graphs for a particular given 
metric dimension is an NP-hard problem. 
This concept has applications in many different domains, 
including graph theory, network architecture, and facility location problems. 
A graph $G$ with order $n$ is known as a Toeplitz graph over the subset $S$ 
of consecutive collections of integers from one to $n$, and two vertices 
will be adjacent to each other if their absolute difference 
is a member of $S$. A graph $G(\mathbb{Z}_{n})$ 
is called a zero-divisor graph over the zero divisors
of a commutative ring $\mathbb{Z}_{n}$, in which two vertices 
will be adjacent to each other if their product 
will leave the remainder zero under modulo $n$. 
Since the local fractional metric dimension problem 
is NP-hard, it is computationally difficult to identify 
an optimal solution or to precisely determine 
the minimal size of a local resolving set; in the worst case, 
the process takes exponential time. Different upper bound 
sequences of local fractional metric dimension are suggested in this article, 
along with a comparison analysis for certain families of Toeplitz and zero-divisor graphs. 
Furthermore, we note that the analyzed local fractional metric dimension upper bounds 
fall into three metric families: constant, limited, and unbounded.}



\begin{document}


\section{Introduction}
\label{sec:01}

Graphs are mathematical structures made up of vertices (nodes) 
and edges (connections) that are used to illustrate various relationships 
and links between objects or entities. 
The need to understand the structural qualities of graphs in terms of the labeling of their vertices 
and the connections between their nodes led to the development of the idea of metric dimension. 
Metric dimension has practical implications in network design, graph theory, 
and optimization~\cite{MR4578575,MR4601162,MR3612583}.

By determining the minimum number of landmarks or reference points required to uniquely 
identify locations in a network, we can efficiently design routing protocols, establish 
communication schemes, and ensure a reliable information flow. Metric dimension also finds 
applications in facility location problems, where the goal is to identify the optimal locations 
for facilities to serve a given set of clients or demands. When one takes into account 
the metric dimension of a graph, we are able to strategically place our facilities 
so that we can serve all of our customers while keeping the overall cost and the amount 
of distance traveled to a minimum. The study of metric dimensions, in general, gives us 
the ability to measure and evaluate the spatial and structural interactions in graphs, 
which in turn leads to practical applications in network design, facility localization, 
sensor networks, graph theory, and the analysis of algorithmic complexity.

Melter and Harary were the ones who originally suggested the idea of metric dimension~\cite{1}. 
Within a network, vertices are denoted by nodes, and the connections that exist between 
the nodes are denoted by edges. This structure makes it possible for agents to travel 
from one vertex to another. To facilitate agent localization, certain points are designated 
as landmarks. The metric basis refers to the smallest set of landmarks, while the metric dimension 
represents the cardinality of this set~\cite{2}. Building upon this foundation, 
Arumugam and Mathew have made significant advancements in the study of the fractional metric 
dimension (FMD)~\cite{3}. Because of their efforts, our basic knowledge of FMD and its 
characteristics has significantly advanced. In addition, the use of FMD has been 
investigated in a variety of network topologies, such as lexicographic, hierarchical, 
Cartesian, corona, and corona-comb networks. These networks are explored utilizing FMD 
collected from product operations, which reveals helpful insights into the structural 
properties of the networks themselves~\cite{4}.

The Fractional Metric Dimension (FMD) in modified Jahangir networks and permutational 
systems was researched by Wang and Feng~\cite{5,6}, as well as by Liu et al.~\cite{7}. 
They investigated the FMD's characteristics as well as various estimate methods while 
applying them to these network architectures. In addition, Raza et al. contributed further 
to the knowledge of FMD across multiple graph classes, while Aisyah et al. calculated 
the boundaries of FMD in metal-organic graphs~\cite{8,9}. 
They recommended FMD to be used for the corona product of systems and proposed 
a more localized version of the metric dimension 
that they called the Local Fractional Metric Dimension (LFMD). 
Their research primarily focused on investigating the LFMD and its applications 
in a variety of graph compositions and formats. The researchers Liu et al. 
made significant contributions to the field of LFMD~\cite{10}. 
More recent work was conducted by Ali, Falc\'{o}n, and other researchers, 
who studied the LFMD of a collection of circular symmetric planar graphs. 
This research was carried out in order to investigate further the local 
fractional metric dimension~\cite{11}. The graphs arise as a result 
of the merging of the edge of an order $n$ cycle with the edges of $n$ 
additional chorded cycles of the same order $n$.
For more on the subject, see~\cite{MR4225551} and references therein.

The research carried out by Aisyah et al. contains one of the first known references to LFMD  
within the context of the corona product of systems~\cite{9}. Within the corona product of graphs, 
the localized fractional metric dimension is taken into account by the LFMD. As a result, 
it is possible to obtain insights on the localized structure and labeling of vertices 
in this graph composition~\cite{13}. The Corona technique, performed in 2016
by Rodr\'{\i}guez-Vel\'{a}zquez et al., introduces the idea of using 
local metric dimensions in graphs~\cite{14}. In the same year, Marsidi et al. 
established the local metric dimensions for line graphs, as well as other different graphs~\cite{15}. 

If all the networks that make up a family, represented by $F$, have the same metric dimension 
within a collection of linked networks, then it is said that the family has a constant metric~\cite{F1}. 
Because a generalized Peterson anti-prism in a $P (n, 2)$ network $C$ is a circulant network $C_2n$, 
each one of them also forms a family with a constant metric dimension. In the course of the inquiry, 
we also look at the size of wheels and Jahangir networks using metric systems and, 
with the help of Integer Programming, we acquire a response that is more specific to our idea. 

FMD has a significant number of innovative features. For instance,
FMD was also employed by Fehr et al. to produce a better solution to a linear programming 
relaxation~\cite{F2}. Afterward, Arumugam and Mathew's contributions helped 
to shed light on the subject~\cite{3}. Moreover, several graph theorists reported a number 
of findings on the metric dimension and graph labeling on many families 
of graphs~\cite{Y2,11,Y6}. Other researchers have studied FMD 
with respect to various other elements of graph theory 
\cite{j1,j2,j3,j4,j5}.

\begin{figure}[H]
\includegraphics[scale=0.7]{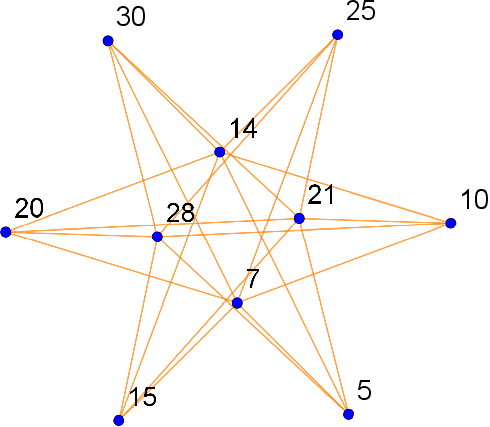}
\caption{The zero-divisor graph for $n=35$.}
\label{w1}
\end{figure}

Algebraic graph theory is a branch of mathematics 
that examines and evaluates graph characteristics and structures 
using algebraic methods. In the field of algebraic graph theory, various problems are still open. For example, the problem of the number of components of a graph, which depends on the modular relation, 
is still dealt with as a conjuncture. There are several algebraic graphs based on the algebraic 
structure that have now been studied, and we refer the reader to~\cite{MR4595248,MR4571844,MR4477215} 
and references therein. Here, we discuss a zero-divisor graph that depends on the set of zero divisors 
of a ring $R$. A graph $\widehat{G}$ is known as a zero-divisor graph whose vertex set is the zero divisors 
of modular ring $\mathbb{Z}_{n}$, and two vertices will be adjacent to each other if their product 
will be zero under $\mod n$~\cite{27N}. The zero-divisor graph for $n=35$ is shown in Figure~\ref{w1}.
\vspace{-6pt}

\begin{figure}[H]
\includegraphics[scale=0.7]{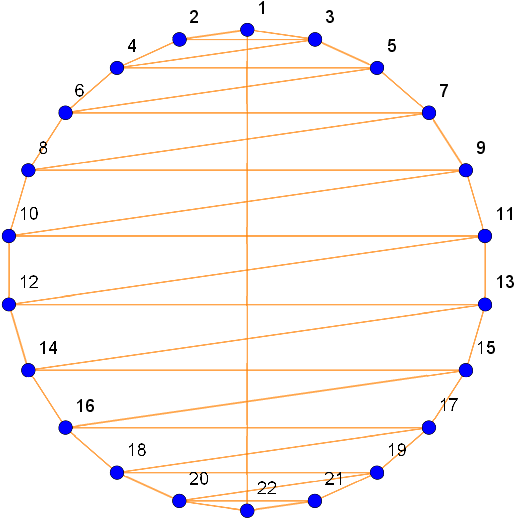}\\
\caption{The Toeplitz graph $T_{22}<\{1, 2, 21\}>$.}\label{fig2}
\end{figure}

A graph $\widetilde{G}$ is called a Toeplitz graph if its vertex set is $M=\{1, 2, \ldots, n\}$ 
and two vertices $p,~q$ are adjacent to each other if their absolute difference belongs to $S$, 
where $S \subseteq M$ ($|p-q|\in S\subseteq M$). The Toeplitz graph $T_{22}<S>$ for $S=\{1, 2, 21\}$ 
is shown in Figure~\ref{fig2}.

A resolving function~\cite{3} 
of the graph $G$ is any map $\zeta: V(G) \rightarrow [0, 1]$ such that
\begin{equation}
\label{eq_RF} 
\sum_{u \in \mathcal{R}\{v, w\} } \zeta(u)\geq 1,
\end{equation}
for every pair of distinct vertices $v,w\in V(G)$. 
The fractional metric dimension of the graph $G$ is
\[
\dim_{\mathrm{f}}(G):= \min\left\{ \sum_{v \in V(G)} \zeta(v) 
\colon\, \zeta \textrm{ is a resolving function of } G \right\}.
\]
The concepts of local resolving neighborhood and local resolving function arise 
similarly, in the case of dealing with pairs of adjacent vertices only. Then, 
the local fractional metric dimension of the graph $G$ is defined as
\[
\mathrm{ldim}_{\mathrm{f}}(G):= \min\left\{ \sum_{v \in V(G)} 
\zeta(v) \colon\, \zeta \textrm{ is a local resolving function of } G \right\}.
\]
Further, from now on, we denote
$\ell(G)=\min\{\left|\mathcal{R}\left\{v,w\right\}\right|\colon\, vw\in E(G)\}$. 
In particular, since $v,w\in\mathcal{R}\{v,w\}$, for all $v,w\in V(G)$, 
one has $\ell(G)\geq 2$. The next result follows from all the previous definitions.
Lemma~\ref{lemma} and Theorem \ref{NT} have a vital importance throughout the manuscript.
In particular, it will be used in the proofs of Theorems~\ref{mt:01} to \ref{thm:07}.

\begin{Lemma}[\cite{9,222}]
\label{lemma} 
If $\widehat{G}$ is any finite simple graph with order greater than $n>2$, then
\begin{equation}
\label{rel}
\frac n{n-\mathrm{ldim}(\widehat{G})+1}
\leq \mathrm{ldim}_{\mathrm{f}}(\widehat{G})
\leq\frac n {\ell(\widehat{G})},
\end{equation}
and the local fractional metric dimension is one if, and only if, 
the graph $\widehat{G}$ is a bipartite graph.
\end{Lemma}

\begin{Theorem}[\cite{LN1}]
\label{NT}
Let $\widehat{G}$ be a simple graph with order $n$ 
and $\mathfrak{R}\{e\}$ be the local resolving neighborhood set. Then,
$$ 
\frac{n}{\beta(\widehat{G})}\leq\mathrm{ldim}_{\mathrm{f}}(\widehat{G}), 
$$
where $\beta(\widehat{G})=\max\{\left|\mathcal{R}\left\{e\right\}\right|
\colon\, e\in E(\widehat{G})\},~~2\leq \beta(\widehat{G})\leq n$.
\end{Theorem}

\begin{figure}[H]
\includegraphics[scale=0.7]{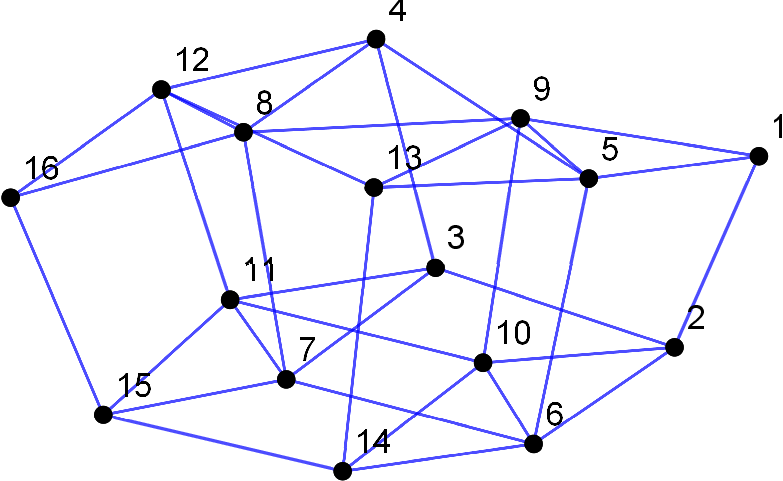} 
\caption{The Toeplitz graph $T_{16}<\{1, 2, 8\}>$.}
\label{Fig1}
\end{figure}

Let us illustrate relation \eqref{rel}.
The Toeplitz graph $G=T_{16}<\{1, 4, 8\}>$ is given in Figure~\ref{Fig1}.
The resolving sets $\mathfrak{R}$, crossposting 
to each pair of adjacent vertices, are given as follows:
\begin{align*}
&\mathfrak{R}\{1,2\}=V\{G\}\setminus\{4,12\}, 
&&\mathfrak{R}\{1,5\}=V\{G\} \backslash \{3,8,9,10,11,16\}, \\
&\mathfrak{R}\{2,3\}=V\{G\}\setminus\{5,13\}, 
&&\mathfrak{R}\{2,6\}=V\{G\}\backslash\{4,9,10,11,12\},\\
&\mathfrak{R}\{3,4\}=V\{G\}\setminus\{1,6,14\}, 
&&\mathfrak{R}\{3,7\}=V\{G\}\backslash\{5,10,11,12,13\},\\
&\mathfrak{R}\{4,5\}=V\{G\}\setminus\{2,7,15\}, 
&&\mathfrak{R}\{4,8\}=V\{G\}\backslash\{1,6,11,12,13,14\},\\
\end{align*}
\begin{align*}
&\mathfrak{R}\{5,6\}=V\{G\}\setminus\{3,8,16\}, 
&&\mathfrak{R}\{5,9\}=V\{G\}\backslash\{1,2,7,12,13,14,15\},\\
&\mathfrak{R}\{6,7\}=V\{G\}\setminus\{4,9\}, 
&&\mathfrak{R}\{6,10\}=V\{G\}\backslash\{1,2,3,8,13,14,15,16\},\\
&\mathfrak{R}\{7,8\}=V\{G\}\setminus\{5,10\}, 
&&\mathfrak{R}\{7,11\}=V\{G\}\backslash\{1,2,3,4,9,14,15,16\},\\
&\mathfrak{R}\{8,9\}=V\{G\}\setminus\{6,11\}, 
&&\mathfrak{R}\{8,12\}=V\{G\}\backslash\{2,3,4,5,10,15,16\},\\
&\mathfrak{R}\{9,10\}=V\{G\}\setminus\{7,12\}, 
&&\mathfrak{R}\{9,13\}=V\{G\}\backslash\{3,4,5,6,11,16\},\\
&\mathfrak{R}\{10,11\}=V\{G\}\setminus\{8,13\}, 
&&\mathfrak{R}\{10,14\}=V\{G\}\backslash\{4,5,6,7,12\},\\
&\mathfrak{R}\{11,12\}=V\{G\}\setminus\{1,9,14\}, 
&&\mathfrak{R}\{11,15\}=V\{G\}\backslash\{5,6,7,8,13\},\\
&\mathfrak{R}\{12,13\}=V\{G\}\setminus\{2,10,15\}, 
&&\mathfrak{R}\{12,16\}=V\{G\}\backslash\{1,6,7,8,9,14\},\\
&\mathfrak{R}\{13,14\}=V\{G\}\setminus\{3,11,16\}, 
&&\mathfrak{R}\{4,12\}=V\{G\}\backslash\{7,8,9\},\\
&\mathfrak{R}\{14,15\}=V\{G\}\setminus\{4, 12\},
&& \mathfrak{R}\{5,13\}=V\{G\}\backslash\{8,9,10\},\\
&\mathfrak{R}\{15,16\}=V\{G\}\setminus\{5, 13\}, 
&& \mathfrak{R}\{6,14\}=V\{G\}\backslash\{9,10,11\},\\
&\mathfrak{R}\{1,9\}=V\{G\}\backslash\{4,5,6\}, 
&&\mathfrak{R}\{7,15\}=V\{G\}\backslash\{10,11,12\},\\
&\mathfrak{R}\{2,10\}=V\{G\}\backslash\{5,6,7\}, 
&& \mathfrak{R}\{8,16\}=V\{G\}\backslash\{11,12,13\},\\
&\mathfrak{R}\{3,11\}=V\{G\}\backslash\{6,7,8\}.
\end{align*}
The cardinalates of the above resolving set are
\begingroup\makeatletter\def\f@size{9}\check@mathfonts
\def\maketag@@@#1{\hbox{\m@th\normalsize \normalfont#1}}%
\begin{align*}
& |\mathfrak{R}\{1,2\}|=14,   &&|\mathfrak{R}\{8,9\}|=14,   
&&|\mathfrak{R}\{15,16\}|=14, &&|\mathfrak{R}\{7,11\}|=8,  &&|\mathfrak{R}\{2,10\}|=13,\\
&|\mathfrak{R}\{2,3\}|=14, &&|\mathfrak{R}\{9,10\}|=14,  
&&|\mathfrak{R}\{1,5\}|=10,   &&|\mathfrak{R}\{8,12\}|=9,  &&|\mathfrak{R}\{3,11\}|=13,\\
&|\mathfrak{R}\{3,4\}|=13, &&|\mathfrak{R}\{10,11\}|=14, 
&&|\mathfrak{R}\{2,6\}|=11,   &&|\mathfrak{R}\{9,13\}|=10,  &&|\mathfrak{R}\{4,12\}|=13,\\
&|\mathfrak{R}\{4,5\}|=13, &&|\mathfrak{R}\{11,12\}|=13, 
&&|\mathfrak{R}\{3,7\}|=11,   &&|\mathfrak{R}\{10,14\}|=11, &&|\mathfrak{R}\{5,13\}|=13,\\
&|\mathfrak{R}\{5,6\}|=13, &&|\mathfrak{R}\{12,13\}|=13, 
&&|\mathfrak{R}\{4,8\}|=10,   &&|\mathfrak{R}\{11,15\}|=11, &&|\mathfrak{R}\{6,14\}|=13,\\
&|\mathfrak{R}\{6,7\}|=14, &&|\mathfrak{R}\{13,14\}|=13, 
&&|\mathfrak{R}\{5,9\}|=9,   &&|\mathfrak{R}\{12,16\}|=10, &&|\mathfrak{R}\{7,15\}|=13,\\
&|\mathfrak{R}\{7,8\}|=14, &&|\mathfrak{R}\{14,15\}|=14, 
&&|\mathfrak{R}\{6,10\}|=8,  &&|\mathfrak{R}\{1,9\}|=13,   &&|\mathfrak{R}\{8,16\}|=13.\\
\end{align*}
\endgroup

Since the minimum and maximum cardinalities of the resolving 
sets are 8 and 14, respectively, it follows by relation \eqref{rel} 
and Theorem~\ref{NT} that
$$
\frac{8}{7}=\frac{16}{14}\leq\mathrm{ldim}_{\mathrm{f}}
\left(T_{16}<1, 4, 8>\right) \leq\frac{16}{8}=2.
$$

The article is arranged as follows. Section~\ref{sec:02} recalls some applications 
of metric dimension in diverse fields such as navigation, networking, pattern recognition 
and image processing, combinatorial optimization, image programming, chemistry, 
and drug discovery. In Section~\ref{sec:03}, stream values of the LFMD of certain families 
of Toeplitz graphs are studied. Then, in Section~\ref{sec:04}, the LFMD for zero-divisor graphs 
over a zero divisor of $\mathds{Z}_{n}$ is computed while in Section~\ref{sec:05}, we compute it
for zero-divisor graphs over $\mathds{Z}_{n}\setminus\{0\}$. In Section~\ref{sec:06}, we give 
some remarks on finding LFMD for algebraic structure graphs. 
In Section~\ref{sec:07}, a few concluding remarks, directions for future work, 
and some open problems needing further research are discussed.


\section{Survey of Applications of Metric Dimension}
\label{sec:02}

Metric dimension, a distance-based characteristic, is used to pinpoint the positions 
of things in space, such as machinery, robots, chemical compounds, etc. The objective 
is to minimize the number of locations or nodes that are used by these items, as well 
as to optimize time consumption and the shortest distances between destinations. 
Applications for metric dimension may be found in the following fields:
\begin{itemize}
\item Navigation, to determine the metric dimension in robotics and navigation 
systems, to set up robots or machines most effectively 
for effective movement and path planning~\cite{MR4603576}.

\item Networking, where metric dimension is used to estimate the bare minimum 
of nodes required to uniquely identify the network~\cite{MR4569822}.

\item Pattern Recognition and Image Processing, 
where metric dimension helps with tasks such as feature 
selection and pattern recognition by allowing data 
points to be represented in a lower-dimensional 
space while still retaining 
their key properties~\cite{MR4225551}.

\item Combinatorial Optimization, to determine the most effective approach 
to cover or represent a collection of points in space, where optimization challenges 
require the usage of metric dimension~\cite{MR4490570}.

\item Integer Programming, where the goal is to discover 
integer solutions that optimize a particular criterion 
\cite{MR4564013}.

\item Chemistry and Drug Discovery, where metric dimension is used 
in molecular networks and chemical compound analysis to characterize and 
describe chemical compounds in a distinctive way, assisting in the identification 
of novel medications or comprehending chemical interactions~\cite{MR4277472}.

\end{itemize}

Overall, metric dimension is a parameter that may be used to solve optimization 
and representation issues in a variety of fields, allowing for effective 
resource allocation and decision-making while taking distance into account. 
It is intriguing to observe how this idea is used in disciplines 
as different as chemistry and computer science.


\section{Stream Values of LFMD for Certain Families of Toeplitz Graphs}
\label{sec:03}

In this section, we find upper bounds sequences of LFMD 
for certain families of Toeplitz graphs.

\begin{Theorem}
\label{mt:01}
Let $T_{\wp}<1, \wp-2>$ be a Toeplitz graph. 
The local fractional metric dimension of $T_{\wp}<1, \wp-2>$ is
\[
\begin{cases}
\begin{array}{ll}
\mathrm{ldim}_{\mathrm{f}}(T_{<1, \wp-2>})=1,
& \text{ if } \wp\in O^{+},\\
\frac{\wp}{\wp-1}\leq\mathrm{ldim}_{\mathrm{f}}(T_{<1, \wp-2>})
\leq\frac{\wp}{\wp-2},&\text{ if } \wp\in E^{+},
\end{array}
\end{cases}
\]
where $O^{+}$ and $E^{+}$ are sets of positive odd integers 
and positive even integers, respectively. 
\end{Theorem}

\begin{proof}
If $\wp$ is an odd number, then the Toeplitz graph's vertex set 
$T_{\wp}<1, \wp-2>$ is divided into two separate sets in the following manner:
\begin{eqnarray*}
S&=&\{1, 3, 5, \ldots, \wp\},  \\
T&=&\{2, 4, 6, \ldots, \wp-1\}.
\end{eqnarray*}
These two sets form a bipartite graph. So, according to Lemma~\ref{lemma}, 
we have 
$$
\mathrm{ldim}_{\mathrm{f}}(T_{<1, \wp-2>})=1.
$$

If $\wp$ is even, then the resolving sets are
$$
\mathfrak{R}\{\vartheta_1,\vartheta_2\}
=\begin{cases}
\begin{array}{ll}
\mathrm{V}(T_{\wp}<1, \wp-2>)-\{{1,\wp}\},
&\text{ if }   \vartheta_1=\frac{\wp}{2}, \vartheta_2=\vartheta_1+1,\\
\mathrm{V}(T_{\wp}<1, \wp-2>)-\{{\frac{\vartheta_1+\vartheta_2}{2}}\},
&\text{ if } \vartheta_1, \vartheta_2 \text{ are both even},\\
\mathrm{V}(T_{\wp}<1, \wp-2>)-\{{\frac{\vartheta_1+\vartheta_2}{2}}\}, 
& \text{ if } \vartheta_1, \vartheta_2 \text{ are both odd},\\
\mathrm{V}(T_{\wp}<1, \wp-2>)-\{{\frac{\wp}{2}+\vartheta_1}\}, 
& \text{ if } \vartheta_2=\vartheta_1+1, \vartheta_1<\frac{\wp}{2},\\
\mathrm{V}(T_{\wp}<1, \wp-2>)-\{{\vartheta_2-\frac{\wp}{2}}\}, 
& \text{ if } \vartheta_2=\vartheta_1+1, \vartheta_1>\frac{\wp}{2}.
\end{array}
\end{cases}
$$
These resolving sets have the following cardinalities:
\begin{center}
$|\mathfrak{R}\{\vartheta_1,\vartheta_2\}|=\begin{cases}\begin{array}{ll}
\wp-2,&\text{ if } \vartheta_1=\frac{\wp}{2}, \vartheta_2=\vartheta_1+1,\\
\wp-1,&\text{ if } \vartheta_1, \vartheta_2 \text{ are both even},\\
\wp-1, & \text{ if } \vartheta_1, \vartheta_2 \text{ are both odd},\\
\wp-1& \text{ if } \vartheta_2=\vartheta_1+1, \vartheta_1<\frac{\wp}{2},\\
\wp-1, & \text{ if } \vartheta_2=\vartheta_1+1, \vartheta_1>\frac{\wp}{2}.
\end{array}\end{cases}$\\
\end{center}
Since the minimum and maximum cardinalities of the 
resolving sets are $\wp-2$ and $\wp-1$, respectively, thus, 
by relation \eqref{rel} and Theorem~\ref{NT}, we have
$$
\frac{\wp}{\wp-1}\leq\mathrm{ldim}_{\mathrm{f}}(T_{<1, \wp-2>})
\leq\frac{\wp}{\wp-2}.
$$
The proof is complete.
\end{proof}

\begin{Theorem}
Let $G=T_{\wp}<1, 2, \wp-1>$ be a Toeplitz graph. 
Then, the local fractional metric dimension is given by 
\[
\begin{cases}
\begin{array}{ll}
\mathrm{ldim}_{\mathrm{f}}(T_\wp<1, 2, \wp-1>)=2,
& \text{ if } \wp=4,\\
\frac{\wp}{\wp-2}\leq{ldim}_{f}(T_\wp<1, 2, \wp-1>)\leq 2,
&\text{ if }  \wp\equiv 0 \mod 4, \wp \neq 4,\\
\frac{\wp}{\wp-1}\leq{ldim}_{f}(T_\wp<1, 2, \wp-1>)\leq \frac{2\wp}{\wp+1},
&\text{ if }  \wp\equiv 1 \mod 4,\\
1\leq{ldim}_{f}(T_\wp<1, 2, \wp-1>)\leq 2,
&\text{ if }  \wp\equiv 2 \mod 4,\\
\frac{\wp}{\wp-1}\leq{ldim}_{f}(T_\wp<1, 2, \wp-1>)\leq \frac{4\wp}{3\wp-5},
&\text{ if }  \wp\equiv 3 \mod 4.\\
\end{array}
\end{cases}
\]
\end{Theorem}

\begin{proof}
There are five cases to be studied.

Case 1. When the Toeplitz graph $G=T_\wp<1, 2, \wp-1>$ has order $\wp=4$, then
\begin{itemize}
\item If  $|\vartheta_1-\vartheta_2|=1,$ then the resolving sets are
\end{itemize}
\[
\mathrm{ }\mathfrak{R}\{\vartheta_1,\vartheta_1+1\} 
= \begin{cases}\begin{array}{ll}
\mathrm{V}(G)\setminus\{3,4\}, & \text{ if } \vartheta_1=1 , \\
\mathrm{V}(G)\setminus\{1,4\}, & \text{ if } \vartheta_1=2 ,\\
\mathrm{V}(G)\setminus\{1,2\}, & \text{ if } \vartheta_1=3 . \\
\end{array}
\end{cases}
\]
\begin{itemize}
\item If  $|\vartheta_1-\vartheta_2|=2,$ then the resolving sets are
\end{itemize}
\[
\mathfrak{R}\{\vartheta_1,\vartheta_1+2\} 
= 
\begin{cases}
\begin{array}{ll}
\mathrm{V}(G)\setminus\{2,4\}, & \text{ if } \vartheta_1=1 ,\\
\mathrm{V}(G)\setminus\{1,3\}, & \text{ if } \vartheta_1=2. \\
\end{array}\end{cases}\]
\begin{itemize}
\item If  $|\vartheta_1-\vartheta_2|=\wp-1,$ then the resolving sets are
\end{itemize}
\[
\mathrm{ }\mathfrak{R}\{\vartheta_1,\wp\} 
= 
\begin{cases}
\begin{array}{ll}
\mathrm{V}(G)\setminus\{2,3\}, 
& \text{ if } \vartheta_1=1. \\
\end{array}
\end{cases}
\]
The cardinalities  of the above resolving set are
\begin{itemize}
\item If  $|\vartheta_1-\vartheta_2|=1,$ then cardinalities of resolving sets are
\end{itemize}
\[
\mathrm{ }|\mathfrak{R}\{\vartheta_1,\vartheta_1+1\}| 
= \begin{cases}
\begin{array}{ll}
\wp-2, & \text{ if } \vartheta_1=1 , \\
\wp-2, & \text{ if } \vartheta_1=2 ,\\
\wp-2, & \text{ if } \vartheta_1=3 . \\
\end{array}
\end{cases}
\]
\begin{itemize}
\item If  $|\vartheta_1-\vartheta_2|=2,$ then cardinalities of resolving sets are
\end{itemize}
\[
|\mathfrak{R}\{\vartheta_1,\vartheta_1+2\}| 
= 
\begin{cases}
\begin{array}{ll}
\wp-2, & \text{ if } \vartheta_1=1 ,\\
\wp-2, & \text{ if } \vartheta_1=2. \\
\end{array}
\end{cases}
\]
\begin{itemize}
\item If $|\vartheta_1-\vartheta_2|=\wp-1$, 
then cardinalities of resolving sets are
\end{itemize}
\[
|\mathfrak{R}\{\vartheta_1,\wp\}| 
= \begin{cases}
\begin{array}{ll}
\wp-2, & \text{ if } \vartheta_1=1.
\end{array}
\end{cases}
\]
Since the minimum cardinality is $\wp-2$,
by relation (\ref{rel}), we have 
$$
{ldim}_{f}\left(T_\wp<1, 2, \wp-1>\right)
\leq \frac{\wp}{\wp-2}=\frac{4}{2}=2.
$$

Case 2. When the Toeplitz graph $T_\wp<1, 2, \wp-1>$ 
has order $\wp\equiv 0 \mod 4, \wp \neq 4$, then
\begin{itemize}
\item If  $|\vartheta_1-\vartheta_2|=1,$ the resolving sets are
\end{itemize}
$\mathfrak{R}\{\vartheta_1,\vartheta_1+1\} =$\\
\begin{adjustwidth}{-\extralength}{0cm}
\[
\begin{cases}
\begin{array}{ll}
\mathrm{V}(G)\setminus\bigcup^{\frac{\wp}{4}}_{\vartheta_2=1}\{2{\vartheta_2}+1\}
\bigcup\{\frac{\wp+4}{2}\}, & \text{ if } \vartheta_1=1,   \\
\mathrm{V}(G)\setminus\{\bigcup^{\frac{\vartheta_1}{2}}_{\vartheta_2=1}\{2{\vartheta_2}-1\}
\bigcup^{\frac{\wp}{4}}_{\vartheta_2=\frac{\vartheta_1+2}{2}}
\{2{\vartheta_2}\}\bigcup^{\frac{\vartheta_1-2}{2}}_{\vartheta_2=0}\\
\{\frac{\wp+4}{2}+2{\vartheta_1}\}\bigcup^{\frac{\wp
-(2\vartheta_1+4)}{2}}_{\vartheta_2=0}\{\wp-\vartheta_2\},
& \text{ if } $1$<{\vartheta_1}\leq{\frac{\wp-2}{2}} , \vartheta_1 \equiv 0 \mod 2 ,\\
\mathrm{V}(G)\setminus\bigcup^{\frac{\vartheta_1-1}{2}}_{\vartheta_2=1}\{2{\vartheta_2}\}
\bigcup^{\frac{\wp}{4}}_{\vartheta_2=\frac{\vartheta_1+1}{2}}\{2{\vartheta_2}+1\}\\
\bigcup^{\frac{\vartheta_1-3}{2}}_{\vartheta_2=0}\{\frac{\wp+6}{2}
+2{\vartheta_2}\}  \bigcup\{\frac{\wp+8}{2}+\vartheta_1-3\}, 
& \text{ if }  1<{\vartheta_1}\leq{\frac{\wp-2}{2}},\vartheta_1 \equiv 1 \mod 2,\\
\mathrm{V}(G)\setminus\bigcup^{\frac{\wp}{4}}_{\vartheta_2=1}\{2{\vartheta_2}-1\}
\bigcup^{\frac{\wp}{4}-1}_{\vartheta_2=0}\{\wp-2{\vartheta_2}\},
&\text{ if } \vartheta_1={\frac{\wp}{2}},\\
\mathrm{V}(G)\setminus\bigcup^{\frac{\vartheta_1}{2}}_{\vartheta_2=1}\{2{\vartheta_2}-1\}
\bigcup^{\frac{\wp}{2}}_{\vartheta_2=\frac{\vartheta_1+2}{2}}\{2{\vartheta_2}\},
& \text{ if }  {\frac{\wp}{2}+1}\leq{\vartheta_1}\leq{\wp-1}, \vartheta_1 \equiv 0 \mod 2  \\
\mathrm{V}(G)\setminus\{\vartheta_1-\frac{\wp}{2}\}\bigcup^{\frac{2\vartheta_1-2}{4}}_{\vartheta_2
=\frac{2\vartheta_1-\wp+2}{4}}\{2{\vartheta_2}\}\bigcup^{
\frac{\wp-2}{2}}_{\vartheta_2=\frac{\vartheta_1+1}{2}}\{2\vartheta_2+1\}, 
& \text{ if } {\frac{\wp}{2}+1}\leq{\vartheta_1}\leq{\wp-2}, \vartheta_1 \equiv 1 \mod 2, \\
\mathrm{V}(G)\setminus\{\vartheta_1-\frac{\wp}{2}\}\bigcup^{\frac{2\vartheta_1-2}{4}}_{\vartheta_2
=\frac{2\vartheta_1-\wp+2}{4}}\{2{\vartheta_2}\}, & \text{ if }  \vartheta_1=\wp-1.
\end{array}
\end{cases}
\]
\end{adjustwidth}
\begin{itemize}
\item If  $|\vartheta_1-\vartheta_2|=2,$ the resolving sets are
\end{itemize}
\begin{adjustwidth}{-\extralength}{0cm}
\[
\mathrm{ }\mathfrak{R}\{\vartheta_1,\vartheta_1+2\}\ 
= 
\begin{cases}
\begin{array}{ll}
\mathrm{V}(G)\setminus\{2\vartheta_1\} \bigcup\{\frac{\wp+4}{2}\}\bigcup\{\frac{\wp+6}{2}\},
& \text{ if } \vartheta_1=1,\\
\mathrm{V}(G)\setminus\{\vartheta_1+1\}\bigcup\{\frac{\wp+2\vartheta_1+4}{2}\},
& \text{ if }  1<{\vartheta_1}<\frac{\wp-2}{2} , \vartheta_1 \equiv 0 \mod 2, \\
\mathrm{V}(G)\setminus\{\vartheta_1+1\}\bigcup\{\frac{\wp+2\vartheta_1+2}{2}\}\bigcup\{
\frac{\wp+2\vartheta_1+4}{2}\}, & \text{ if }  1<{\vartheta_1}
<\frac{\wp-2}{2} ,  \vartheta_1 \equiv 1 \mod 2, \\
\mathrm{V}(G)\setminus\{\vartheta_1+1\}\bigcup\{\frac{\wp+2\vartheta_1+2}{2}\},
& \text{ if } \vartheta_1=\frac{\wp-2}{2}, \\
\mathrm{V}(G)\setminus\{\frac{\wp}{2\vartheta_1}\}\bigcup\{\vartheta_1+1\},
& \text{ if }  \vartheta_1=\frac{\wp}{2}, \\
\mathrm{V}(G)\setminus\{\vartheta_1-\frac{\wp-1}{2}\}\bigcup\{\vartheta_1+1\},
& \text{ if } \frac{\wp}{2}<{\vartheta_1}<\wp-2, \vartheta_1 \equiv 1 \mod 2, \\
\mathrm{V}(G)\setminus\{\frac{2\vartheta_1-\wp}{2}\}
\bigcup\{\frac{2\vartheta_1-\wp+2}{2}\}\bigcup\{\vartheta_1+1\}, 
& \text{ if }  \frac{\wp}{2}<{\vartheta_1}<\wp-2, 
\vartheta_1 \equiv 0 \mod 2, \\
\mathrm{V}(G)\setminus\{\frac{2\vartheta_1-\wp}{2}\}
\bigcup\{\frac{\vartheta_1}{2}\}\bigcup\{\vartheta_1+1\}, 
& \text{ if } {\vartheta_1}\equiv 1 \mod 2.
\end{array}
\end{cases}
\]
\end{adjustwidth}
\begin{itemize}
\item If $|\vartheta_1-\vartheta_2|=\wp-1$, the resolving sets are
\end{itemize}
$$
\mathfrak{R}\{1,\wp\}=\left\{\frac{\wp}{2}\right\}
\bigcup\left\{\frac{\wp+2}{2}\right\}.
$$
The cardinalities of the above resolving set are
\begin{itemize}
\item If  $|\vartheta_1-\vartheta_2|=1,$ then
\end{itemize}
\[
|\mathfrak{R}\{\vartheta_1,\vartheta_1+1\}| 
= 
\begin{cases}
\begin{array}{ll}
\frac{3\wp-4}{4}, & \text{ if } \vartheta_1=1,   \\
\frac{\wp+4\vartheta_1}{4}, & \text{ if }  $1$<{\vartheta_1}\leq{\frac{\wp-2}{2}},\\
\frac{3\wp-2\vartheta_1-4}{4},& \text{ if }  1<{\vartheta_1}\leq{\frac{\wp-2}{2}}, \\
\frac{\wp}{2},&\text{ if } \vartheta_1={\frac{\wp}{2}},\\
\frac{3\wp-4}{4},& \text{ if }  {\frac{\wp}{2}+1}\leq{\vartheta_1}\leq{\wp-1}, \\
\frac{\wp+2\vartheta_1-2}{4} ,& \text{ if } {\frac{\wp}{2}+1}\leq{\vartheta_1}\leq{\wp-2}, \\
\frac{3\wp-4}{4},& \text{ if }  \vartheta_1={\wp-1}. 
\end{array}
\end{cases}
\]
\begin{itemize}
\item If  $|\vartheta_1-\vartheta_2|=2,$ then
\end{itemize}
\[
|\mathfrak{R}\{\vartheta_1,\vartheta_1+2\}|\ 
= \begin{cases}
\begin{array}{ll}
\wp-3,& \text{ if } \vartheta_1=1,\\
\wp-2 ,& \text{ if }  1<{\vartheta_1}<\frac{\wp-2}{2} , \vartheta_1 \equiv 0 \mod 2, \\
\wp-3,& \text{ if }  1<{\vartheta_1}<\frac{\wp-2}{2} ,  \vartheta_1 \equiv 1 \mod 2, \\
\wp-2,& \text{ if } \vartheta_1=\frac{\wp-2}{2}, \\
\wp-2,& \text{ if }  \vartheta_1=\frac{\wp}{2} , \\
\wp-2,& \text{ if } \frac{\wp}{2}<\vartheta_1<\wp-2 ,\vartheta_1 \equiv 1 \mod 2, \\
\wp-3, & \text{ if } \frac{\wp}{2}<\vartheta_1<\wp-2 , \vartheta_1 \equiv 0 \mod 2, \\
\wp-3, & \text{ if } {\vartheta_1}\equiv 1 \mod 2. 
\end{array}
\end{cases}
\] 
\begin{itemize}
\item If  $|\vartheta_1-\vartheta_2|=\wp-1,$ then
\end{itemize}
$$
|\mathfrak{R}\{\vartheta_1,\wp\}|= \wp-2.
$$
Since the minimum and maximum cardinalities of the resolving 
sets are $\frac{\wp}{2}$ and $\wp-2$, respectively, 
by relation \eqref{rel} and Theorem~\ref{NT}, we have:
$$
\frac{\wp}{\wp-2}\leq{ldim}_{f}(T_\wp<1, 2, \wp-1>)\leq \frac{\wp}{\wp/2}=2.
$$
Case 3. When the Toeplitz graph $G=T_\wp<1, 2, \wp-1>$ 
has order $\wp\equiv 1 \mod 4$, then the following resolving sets are
\begin{itemize}
\item If  $|\vartheta_1-\vartheta_2|=1$, then 
$\mathfrak{R}\{\vartheta_1,\vartheta_1+1\} =$
\end{itemize}
\begin{adjustwidth}{-\extralength}{0cm}
\[
\begin{cases}
\begin{array}{ll}
\mathrm{V}(G)\setminus\bigcup^{\frac{\wp-1}{4}}_{\vartheta_2=1}\{2{\vartheta_2}+1\},
& \text{ if } \vartheta_1=1,\\
\mathrm{V}(G)\setminus\bigcup^{\frac{\vartheta_1}{2}}_{\vartheta_2=1}\{2{\vartheta_2}-1\}
\bigcup^{\frac{\wp+3}{4}}_{\vartheta_2=\frac{\vartheta_1+2}{2}}\{2{\vartheta_2}\}\\
\bigcup^{\frac{\vartheta_1-2}{2}}_{\vartheta_2=0}\{\frac{\wp+7}{2}
+2{\vartheta_2}\}\bigcup^{\frac{\wp-(2\vartheta_1+1)}{2}}_{\vartheta_2=0}
\{\wp-\vartheta_2\},& \text{ if } 1<{\vartheta_1}\leq{\frac{\wp-1}{2}}, \vartheta_1 \equiv 0 \mod 2,\\
\mathrm{V}(G)\setminus\bigcup^{\frac{\vartheta_1-1}{2}}_{\vartheta_2=1}\{2{\vartheta_2}\}
\bigcup^{\frac{\wp+2\vartheta_1+1}{4}}_{\vartheta_2=\frac{\vartheta_1+3}{2}}
\{2{\vartheta_2}-1\},&\text{ if } 1<{\vartheta_1}\leq{\frac{\wp-1}{2}} , \vartheta_1 \equiv 1 \mod 2,\\
\mathrm{V}(G)\setminus\bigcup^{\frac{\vartheta_1}{2}}_{\vartheta_2=1}\{2{\vartheta_2}-1\}
\bigcup^{\frac{\wp-1}{2}}_{\vartheta_2=\frac{\vartheta_1+2}{2}}\{2{\vartheta_2}\}, 
&\text{ if } \vartheta_1=\frac{\wp-1}{2}, \vartheta_1 \equiv 0 \mod 2,\\
\mathrm{V}(G)\setminus\bigcup^{\frac{\wp-1}{2}}_{\vartheta_2=1}\{2{\vartheta_2}\}
\bigcup^{\vartheta_1}_{\vartheta_2=\vartheta_1+2}\{2{\vartheta_2}-1\}, 
&\text{ if } \vartheta_1 = \frac{\wp+1}{2},\\
\mathrm{V}(G)\setminus\bigcup^{\frac{2\vartheta_1-\wp-1}{4}}_{\vartheta_2=1}\{2{\vartheta_2}-1\}
\bigcup^{\frac{\vartheta_1-1}{2}}_{\vartheta_2=1}\{2{\vartheta_2}\}
\bigcup^{\frac{\wp+1}{2}}_{\vartheta_2=\frac{\wp+3}{2}}\{2{\vartheta_2}-1\}, 
&\text{ if } \frac{\wp+1}{2}<{\vartheta_1}<{\wp-1} , \vartheta_1 \equiv 1 \mod 2,\\
\mathrm{V}(G)\setminus\bigcup^{\frac{\vartheta_1-2}{2}}_{\vartheta_2=\frac{2\vartheta_1-\wp+1}{4}}
\{2{\vartheta_2}+1\}\bigcup^{\frac{\wp-1}{2}}_{\vartheta_2
=\frac{\vartheta_1+2}{2}}\{2{\vartheta_2}\}, 
&\text{ if } \frac{\wp+1}{2}<{\vartheta_1}<{\vartheta_2-1}, 
\vartheta_1 \equiv 0 \mod 2, \\
\mathrm{V}(G)\setminus\bigcup^{\frac{\vartheta_1-2}{2}}_{\vartheta_2
=\frac{\wp-1}{4}}\{2{\vartheta_2}+1\},  
&\text{ if } \vartheta_1=\wp-1. 
\end{array}
\end{cases}
\]
\end{adjustwidth}
\begin{itemize}
\item If $|\vartheta_1-\vartheta_2|=2$, then
\end{itemize}
\[
\mathrm{}\mathfrak{R}\{\vartheta_1,\vartheta_1+2\}
=\begin{cases}
\begin{array}{ll}
\mathrm{V}(G)\setminus\{\vartheta_1+1\}\bigcup\{\frac{\wp+2\vartheta_1+3}{2}\},
& \text{ if } {1}\leq{\vartheta_1}\leq\{\frac{\wp-1}{2}\} , {\vartheta_1}\equiv 1 \mod 2, \\
\mathrm{V}(G)\setminus\{\vartheta_1+1\},&\text{ if } {1}\leq{\vartheta_1}
\leq\{\frac{\wp-1}{2}\} , \vartheta_1 \equiv 0 \mod 2,\\
\mathrm{V}(G)\setminus\{\vartheta_1-\frac{\wp-1}{2}\}\bigcup\{\vartheta_1+1\},
&\text{ if } \frac{\wp-1}{2}<\vartheta_1\leq {\wp-2} , \vartheta_1 \equiv 1 \mod 2,\\
\mathrm{V}(G)\setminus\{\vartheta_1+1\}, 
&\text{ if } \frac{\wp-1}{2}<i
\leq {\wp-2} , \vartheta_1 \equiv 0 \mod 2.
\end{array}
\end{cases}
\]
\begin{itemize}
\item If  $|\vartheta_1-\vartheta_2|=\wp-1,$ then
\end{itemize}
$$
\mathfrak{R}\{1,\wp\}=\left\{\frac{\wp+1}{2}\right\}.
$$ 
The cardinalities of the above resolving sets are
\begin{itemize}
\item If  $|\vartheta_1-\vartheta_2|=1,$ then
\end{itemize}
\[
\mathrm{}|\mathfrak{R}\{\vartheta_1,\vartheta_1+1\}|
= 
\begin{cases}
\begin{array}{ll}
\frac{3\wp-3}{4},& \text{ if } \vartheta_1=1,\\
\frac{\wp+2\vartheta_1-5}{4},
& \text{ if } 1<{\vartheta_1}\leq{\frac{\wp-1}{2}}, \\
\frac{3\wp+2\vartheta_1+10}{4}, 
& \text{ if } 1<{\vartheta_1}\leq{\frac{\wp-1}{2}}, \\
\frac{\wp+1}{2}, 
&\text{ if } \vartheta_1=\frac{\wp-1}{2}, \vartheta_1 \equiv 0 \mod 2,\\
\frac{\wp+3}{2}, 
&\text{ if } \vartheta_1 = \frac{\wp+1}{2},\\
\frac{5\wp-4\vartheta_1+5}{4}, 
&\text{ if } \frac{\wp+1}{2}<\vartheta_1< \wp-1, \vartheta_1 \equiv 1 \mod 2,\\
\frac{3\wp+2\vartheta_1-7}{4}, 
&\text{ if } \frac{\wp+1}{2}<{\vartheta_1}<{\wp-1} , \vartheta_1 \equiv 0 \mod 2,\\
\frac{5\wp-2\vartheta_1-1}{4}, 
&\text{ if }, {\vartheta_1}\leq{\wp-1}. \\
\end{array}
\end{cases}
\] 
\begin{itemize}
\item If  $|\vartheta_1-\vartheta_2|=2,$ then
\end{itemize}
\[
\mathrm{}|\mathfrak{R}\{\vartheta_1,\vartheta_1+2\}|
=\begin{cases}
\begin{array}{ll}
\wp-2,& \text{ if } {1}
\leq{\vartheta_1}\leq\{\frac{\wp-1}{2}\} , {\vartheta_1}\equiv 1 \mod 2, \\
\wp-1,&\text{ if } {1}
\leq{\vartheta_1}\leq\{\frac{\wp-1}{2}\} , \vartheta_1 \equiv 0 \mod 2,\\
\wp-2,&\text{ if } \frac{\wp-1}{2}<\vartheta_1
\leq {\wp-2} , \vartheta_1 \equiv 1 \mod 2,\\
\wp-1, &\text{ if } \frac{\wp-1}{2}<\vartheta_1
\leq {\wp-2} , \vartheta_1 \equiv 0 \mod 2.
\end{array}
\end{cases}
\] 
\begin{itemize}
\item If  $|\vartheta_1-\vartheta_2|=\wp-1,$ then
\end{itemize}
$$
|\mathfrak{R}\{1,\wp\}|=\wp-1.
$$
The minimum and maximum cardinalities of the resolving sets are 
$\frac{\wp+1}{2}$ and $\wp-1$, respectively. Therefore, by relation 
\eqref{rel} and Theorem~\ref{NT}, one has that
$$
\frac{\wp}{\wp-1}\leq{ldim}_{f}(T_\wp<1, 2, \wp-1>)\leq \frac{2\wp}{\wp+1}.
$$
Case 4. When the Toeplitz graph $T_\wp<1, 2, \wp-1>$ 
has order $\wp\equiv 2 \mod 4$, then we have the following resolving sets:
\begin{itemize}
\item If $|\vartheta_1-\vartheta_2|=1$, then
$\mathfrak{R}\{\vartheta_1,\vartheta_1+1\}=$
\end{itemize}
\begin{adjustwidth}{-\extralength}{0cm}
\[
\begin{cases}
\begin{array}{ll}
\mathrm{V}(G)\setminus\bigcup^{\frac{\wp+2}{4}}_{\vartheta_2=1}\{2{\vartheta_2}+1\}, 
& \text{ if } \vartheta_1=1,\\
\mathrm{V}(G)\setminus\bigcup^{\frac{\vartheta_1}{2}}_{\vartheta_2=1}\{2{\vartheta_2}-1\}
\bigcup^{\frac{\wp+2}{4}}_{\vartheta_2=\frac{\vartheta_1+2}{2}}\{2{\vartheta_2}\}\\
\bigcup^{\frac{\vartheta_1-2}{2}}_{\vartheta_2=0}\{\frac{\wp+6}{2}
+2{\vartheta_2}\}\bigcup^{\frac{\wp-(2\vartheta_1+4)}{2}}_{\vartheta_2=0}
\{\wp-\vartheta_2\}, &\text{ if }  1<{\vartheta_1}
<{\frac{\wp}{2}-1} , \vartheta_1 \equiv 0 \mod 2, \\
\mathrm{V}(G)\setminus\bigcup^{\frac{\vartheta_1}{2}}_{\vartheta_2=1}\{2\vartheta_2-1\}
\bigcup^{\frac{\wp+2}{4}}_{\vartheta_2=\frac{\vartheta_1+2}{2}}
\{2\vartheta_2\}\bigcup^{\frac{\vartheta_1-2}{2}}_{\vartheta_2=0}
\{\frac{\wp+6}{2}+2\vartheta_2\},&\text{ if } 1<{\vartheta_1}
\leq\frac{\wp-2}{2},  {\vartheta_1}\equiv 0 \mod 2, \\
\mathrm{V}(G)\setminus\bigcup^{\frac{\vartheta_1-1}{2}}_{\vartheta_2=1}\{2{\vartheta_2} \}
\bigcup^{\frac{\wp+2\vartheta_1}{4}}_{\vartheta_2=\frac{\vartheta_1+1}{2}}
\{2{\vartheta_2}+1\},&\text{ if } ,\vartheta_1
< \frac{\wp}{2}-1,\vartheta_1 \equiv 1 \mod 2, \\
\mathrm{V}(G)\setminus\bigcup^{\frac{\wp-2}{4}}_{\vartheta_2=1}\{2{\vartheta_2}\}
\bigcup^{\frac{\wp-2}{2}}_{\vartheta_2=\frac{\wp+2}{4}}\{2{\vartheta_2}+1\}, 
&\text{ if } \vartheta_1 = \frac{\wp}{2},\\
\mathrm{V}(G)\setminus\bigcup^{\frac{\vartheta_1}{2}}_{\vartheta_2=1}\{2{\vartheta_2}-1\}
\bigcup^{\vartheta_1-1}_{\vartheta_2=\frac{\vartheta_1+2}{2}}\{2{\vartheta_2}\}, 
&\text{ if } {\vartheta_1}={\frac{\wp+2}{2}},\\
\mathrm{V}(G)\setminus\bigcup^{\frac{\vartheta_1}{2}}_{\vartheta_2= 1}\{2{\vartheta_2}-1\}
\bigcup^{\frac{\wp-8}{2}}_{\vartheta_2=1}\{2{\vartheta_2}\}, 
&\text{ if } \frac{\wp+1}{2}<{\vartheta_1}<{\wp-1} , \vartheta_1 \equiv 0 \mod 2, \\
\mathrm{V}(G)\setminus\bigcup^{\frac{\vartheta_1-1}{2}}_{\vartheta_2
=\frac{2\vartheta_1-\wp}{4}}\{2{\vartheta_2}\}\bigcup^{
\frac{\wp-2}{2}}_{\vartheta_2=\frac{\vartheta_1+1}{2}}\{2\vartheta_2+1\}, 
&\text{ if } \frac{\wp+1}{2}<{\vartheta_1}<{\wp-1} , \vartheta_1 \equiv 1 \mod 2, \\
\mathrm{V}(G)\setminus\bigcup^{\frac{\vartheta_1-1}{2}}_{\vartheta_2
=\frac{2\vartheta_1-\wp}{4}}\{2{\vartheta_2}\} , &\text{ if }  \vartheta_1=\wp-1.
\end{array}
\end{cases}
\]
\end{adjustwidth}
\begin{itemize}
\item If  $|\vartheta_1-\vartheta_2|=2,$ then
\end{itemize}
\begin{adjustwidth}{-\extralength}{0cm}
\[
\mathfrak{R}\{\vartheta_1,\vartheta_1+2\}
=
\begin{cases}
\begin{array}{ll}
\mathrm{V}(G)\setminus\{\vartheta_1+1\}, & \text{ if } {1}\leq{\vartheta_1}
\leq\frac{\wp-2}{2} , {\vartheta_1}\equiv 1 \mod 2,\\
\mathrm{V}(G)\setminus\{\vartheta_1+1\}\bigcup\{\frac{\wp+2\vartheta_1+4}{2}\}, 
& \text{ if } {1}\leq{\vartheta_1}\leq\frac{\wp-2}{2} , {\vartheta_1}\equiv 0 \mod 2, \\
\mathrm{V}(G)\setminus\{\vartheta_1,\vartheta_1+2\}=\{\vartheta_1+1\}, 
&\text{ if } \frac{\wp-2}{2}\leq{\vartheta_1}\leq{\frac{\wp+2}{2}},\\
\mathrm{V}(G)\setminus\{\frac{2\vartheta_1-\wp}{2}\}\bigcup\{\vartheta_1+1\},
&\text{ if } \frac{\wp+2}{2}<{\vartheta_1}\leq{\wp-2} , \vartheta_1 \equiv 1 \mod 2, \\
\mathrm{V}(G)\setminus\{\vartheta_1+1\} ,&\text{ if } \frac{\wp+2}{2}
<{\vartheta_1}\leq{\wp-2} , \vartheta_1 \equiv 0 \mod 2,
\end{array}
\end{cases}
\] 
\end{adjustwidth}
for $\wp$ $\equiv 2 \mod 4$.

\begin{itemize}
\item If  $|\vartheta_1-\vartheta_2|=\wp-1$, then we have
$$
\mathfrak{R}\{1,\wp\}=\emptyset.
$$
\end{itemize}

The cardinalities of the above resolving set are 
\begin{itemize}
 \item If  $|\vartheta_1-\vartheta_2|=1,$ then
\end{itemize}
\[
|\mathfrak{R}\{\vartheta_1,\vartheta_1+1\}|
=
\begin{cases}
\begin{array}{ll}
\frac{3\wp-2}{4}, & \text{ if } \vartheta_1=1,\\
\frac{\wp+2\vartheta_1+2}{4}, &\text{ if }  1
<{\vartheta_1}<\frac{\wp}{2}-1,\vartheta_1 \equiv 0 \mod 2, \\
\frac{3\wp-2\vartheta_1-2}{4},&\text{ if } 
1<{\vartheta_1} \leq\frac{\wp}{2}-1, \vartheta_1 \equiv 0 \mod 2, \\
\frac{3\wp-2\vartheta_1}{4}, &\text{ if } \vartheta_1 
<\frac{\wp}{2}-1, \vartheta_1 \equiv 1 \mod 2, \\
\frac{\wp+2}{2}, &\text{ if } \vartheta_1 = \frac{\wp}{2}, \\
\frac{\wp}{2}, &\text{ if } {\vartheta_1}={\frac{\wp+2}{2}}, \\
\frac{2\wp-\vartheta_1+4}{2}, &\text{ if } 
\frac{\wp+1}{2}<{\vartheta_1}<{\wp-1}, \vartheta_1 \equiv 0 \mod 2, \\
\frac{\wp+2e}{4} , &\text{ if }  
\frac{\wp+1}{2}<{\vartheta_1}<{\wp-1} , \vartheta_1 \equiv 1 \mod 2, \\
\frac{3\wp+2}{4} , & \text{ if }  \vartheta_1=\wp-1.
\end{array}
\end{cases}
\] 
\begin{itemize}
\item If  $|\vartheta_1-\vartheta_2|=2,$ then
\end{itemize}
\[
|\mathfrak{R}\{\vartheta_1,\vartheta_1+2\}|
=\begin{cases}
\begin{array}{ll}
\wp-1, & \text{ if }  1\leq \vartheta_1 \leq\frac{\wp-2}{2} , \vartheta_1 \equiv 1 \mod 2,\\
\wp-2, & \text{ if }  1\leq \vartheta_1 \leq\frac{\wp-2}{2} , \vartheta_1 \equiv 0 \mod 2,\\
\wp-1, & \text{ if }  \frac{\wp-2}{2}\leq \vartheta_1 \leq\frac{\wp+2}{2} , \vartheta_1 \equiv 1 \mod 2,\\
\wp-2, & \text{ if }  \frac{\wp+2}{2}< \vartheta_1 \leq \wp-2 , \vartheta_1 \equiv 1 \mod 2,\\
\wp-1, & \text{ if }  \frac{\wp+2}{2}< \vartheta_1 \leq \wp-1 , \vartheta_1 \equiv 0 \mod 2.
\end{array}
\end{cases}
\]
\begin{itemize}
\item If  $|\vartheta_1-\vartheta_2|=\wp-1,$ then
\end{itemize}
$$
|\mathfrak{R}\{1,\wp\}|=\wp.
$$
Since the minimum and maximum cardinalities of the resolving 
sets are $\frac{\wp}{2}$ and $\wp$, respectively, it follows 
by relation \eqref{rel} and Theorem~\ref{NT} that
$$
1=\frac{\wp}{\wp}\leq{ldim}_{f}(T_\wp<1, 2, \wp-1>)
\leq \frac{\wp}{\wp/2}=2.
$$
Case 5. When the Toeplitz graph $G=T_\wp<1, 2, \wp-1>$ 
has order $\wp \equiv 3 \mod 4$, then the resolving sets are
\begin{itemize}
\item If  $|\vartheta_1-\vartheta_2|=1$, 
then $\mathfrak{R}\{\vartheta_1,\vartheta_1+1\}=$
\end{itemize}

\begin{adjustwidth}{-\extralength}{0cm}

\[
\begin{cases}
\begin{array}{ll}
\mathrm{V}(G)\setminus\bigcup^{\frac{\wp+1}{4}}_{\vartheta_2=1}\{2{\vartheta_2}+1\}
\bigcup\{\frac{\wp+5}{2}\},& \text{ if } \vartheta_1=1,\\
\mathrm{V}(G)\setminus\bigcup^{\frac{\vartheta_1}{2}}_{\vartheta_2=1}\{2{\vartheta_2}-1\}
\bigcup^{\frac{\wp+5}{4}}_{\vartheta_2=\frac{\vartheta_1+2}{2}}\{2{\vartheta_2}\}\\
\bigcup^{\frac{\vartheta_1-2}{2}}_{\vartheta_2=0}\{
\frac{\wp+5}{2}+2{\vartheta_2}\}\bigcup^{\frac{\wp
-(2\vartheta_1+3)}{2}}_{\vartheta_2=0}\{\wp-\vartheta_2\},
& \text{ if } 1<{\vartheta_1}\leq{\frac{\wp-3}{2}} , \vartheta_1 \equiv 0 \mod 2, \\
\mathrm{V}(G)\setminus\bigcup^{\frac{\vartheta_1-1}{2}}_{\vartheta_2=1}\{2{\vartheta_2}\}
\bigcup^{\frac{\wp+1}{4}}_{\vartheta_2=\frac{\wp+1}{4}}\{2{\vartheta_2}+1\}\\
\mathrm{V}(G)\setminus\bigcup^{\frac{\vartheta_1-3}{2}}_{\vartheta_2=0}\{\frac{\wp+7}{2}
+2{\vartheta_2}\}\bigcup\{\frac{\wp+9}{2}+(\vartheta_1-3)\},
& \text{ if } 1<{\vartheta_1}\leq{\frac{\wp-3}{2}} , \vartheta_1 \equiv 1 \mod 2, \\
\mathrm{V}(G)\setminus\bigcup^{\frac{\vartheta_1-1}{2}}_{\vartheta_2=1}\{2{\vartheta_2}\}
\bigcup\{\frac{\wp+3}{2}\}\bigcup^{\frac{\wp-3}{4}}_{\vartheta_2=1}
\{\frac{\wp+3}{2}+2\vartheta_2\},
& \text{ if } \frac{\wp-1}{2}\leq{\vartheta_1}\leq\frac{\wp+1}{2},\\
\mathrm{V}(G)\setminus\bigcup^{\frac{\vartheta_1-2}{2}}_{\vartheta_2=1}\{2{\vartheta_2}-1\}
\bigcup\{\frac{\wp-1}{2}\}\bigcup^{\frac{\wp-3}{4}}_{\vartheta_2=1}
\{\frac{\wp+1}{2}+2\vartheta_2\},
& \text{ if } \frac{\wp+1}{2}\leq{\vartheta_1}\leq\frac{\wp+3}{2},\\
\mathrm{V}(G)\setminus\{\vartheta_1-\frac{\wp+1}{2}\}\bigcup^{\frac{2\vartheta_1}{4}}_{\vartheta_2
=\frac{2\vartheta_1-\wp+3}{4}}\{2\vartheta_2-1\}
\bigcup^{\frac{\wp-1}{2}}_{\vartheta_2=\frac{\vartheta_1+2}{2}}\{2\vartheta_2\},
& \text{ if } \frac{\wp}{2}+2<{\vartheta_1}\leq{\wp-2}, \\
\mathrm{V}(G)\setminus\{\vartheta_1-\frac{\wp+1}{2}\}\bigcup^{\frac{2\vartheta_1}{4}}_{\vartheta_2
=\frac{2\vartheta_1-\wp+3}{4}}\{2\vartheta_2-1\},
& \text{ if }   \vartheta_1=\wp-1 , \\
\mathrm{V}(G)\setminus\bigcup^{\frac{2\vartheta_1-\wp+1}{4}}_{\vartheta_2=1}\{2\vartheta_2-1\}
\bigcup^{\frac{\vartheta_1-1}{2}}_{\vartheta_2=1}\{2\vartheta_2\}
\bigcup^{\frac{\wp-1}{2}}_{\vartheta_2=\frac{\vartheta_1+1}{2}}\{2\vartheta_2+1\},
& \text{ if } {\vartheta_1}={\wp-1} , \vartheta_1 \equiv 1 \mod 2.
\end{array}
\end{cases}
\]

\end{adjustwidth}
\begin{itemize}
\item If $|\vartheta_1-\vartheta_2|=2$, then
\end{itemize}
\begin{adjustwidth}{-\extralength}{0cm}
\[
\mathfrak{R}\{\vartheta_1,\vartheta_1+2\}
=
\begin{cases}
\begin{array}{ll}
\mathrm{V}(G)\setminus\{2\vartheta_1\} \bigcup\{\frac{\wp+5}{2}\}, 
& \text{ if }   \vartheta_1=1,\\
\mathrm{V}(G)\setminus\{\vartheta_1+1\} \bigcup\{\frac{\wp+2\vartheta_1+3}{2}\}
\bigcup\{\frac{\wp+2\vartheta_1+5}{4}\}, 
& \text{ if } 1<\vartheta_1<\frac{\wp-3}{2}, 
\vartheta_1 \equiv 0 \mod 2,\\
\mathrm{V}(G)\setminus\{\vartheta_1+1\}\bigcup\{\frac{\wp+2\vartheta_1+3}{2}\},
&\text{ if } \vartheta_1=\frac{\wp-3}{2}, {{\vartheta_1}\equiv 0 \mod 2}, \\
\mathrm{V}(G)\setminus\{\vartheta_1+1\}, &\text{ if } \vartheta_1=\frac{\wp-1}{2}, \\
\mathrm{V}(G)\setminus\{\vartheta_1-\frac{\wp-1}{2}\}\bigcup\{\vartheta_1+1\}, 
&\text{ if } \vartheta_1=\frac{\wp+1}{2} ,\vartheta_1 \equiv 0 \mod 2, \\
\mathrm{V}(G)\setminus\{\frac{2\vartheta_1-\wp-1}{2}\}\bigcup\{\frac{2\vartheta_1
-\wp+1}{2}\}\bigcup\{\vartheta_1+1\}, 
&\text{ if } \frac{\wp+2}{2}<{\vartheta_1}\leq{\wp-3},\vartheta_1 \equiv 0 \mod 2, \\
\mathrm{V}(G)\setminus\{\frac{2\vartheta_1-\wp+1}{2}\}\bigcup\{ \vartheta_1+1\}, 
&\text{ if } \frac{\wp-1}{2}<\vartheta_1<\wp-2 , \\
\mathrm{V}(G)\setminus\{\frac{\wp-1}{2}\}\bigcup\{\vartheta_1+1\},
&\text{ if } \vartheta_1=\wp-2.
\end{array}
\end{cases}
\]
\end{adjustwidth}
\begin{itemize}
\item If $|\vartheta_1-\vartheta_2|=\wp-1,$ then
\end{itemize}
$$
\mathfrak{R}\{1,\wp\}=\left\{\frac{\wp+1}{2}\right\}.
$$
The cardinalities of the above resolving set are
\begin{itemize}
\item If $|\vartheta_1-\vartheta_2|=1$, then
\end{itemize}
\[
|\mathfrak{R}\{\vartheta_1,\vartheta_1+1\}|
=
\begin{cases}
\begin{array}{ll}
\frac{3\wp-5}{4},& \text{ if } \vartheta_1=1,\\
\frac{\wp+2\vartheta_1-3}{4},
& \text{ if } 1<{\vartheta_1}\leq{\frac{\wp-3}{2}}, \vartheta_1 \equiv 0 \mod 2, \\
\wp-1,& \text{ if } 1<{\vartheta_1}\leq{\frac{\wp-3}{2}}, \vartheta_1 \equiv 1 \mod 2, \\
\frac{3\wp-2\vartheta_1+1}{4},& \text{ if } \frac{\wp-1}{2}\leq{\vartheta_1}\leq\frac{\wp+1}{2}, \\
\frac{3\wp-2\vartheta_1+3}{4},& \text{ if } \frac{\wp+1}{2}\leq{\vartheta_1}\leq\frac{\wp+3}{2} , \\
\frac{3\wp+2\vartheta_1-9}{4},& \text{ if } \frac{\wp}{2}+2<{\vartheta_1}\leq{\wp-2}, \\
\frac{3\wp-5}{4},& \text{ if } \vartheta_1=\wp-1, \vartheta_1 \equiv 0 \mod 2, \\
\frac{5\wp-3}{4},& \text{ if } \vartheta_1=\wp-1, \vartheta_1 \equiv 1 \mod 2.
\end{array}
\end{cases}
\] 
\begin{itemize}
\item If  $|\vartheta_1-\vartheta_2|=2$, then
\end{itemize}
\[
|\mathfrak{R}\{\vartheta_1,\vartheta_1+2\}|
=
\begin{cases}
\begin{array}{ll}
\wp-2,& \text{ if } \vartheta_1=1,\\
\wp-3 ,& \text{ if } 1<\vartheta_1<\frac{\wp-3}{2}, \vartheta_1 \equiv 0 \mod 2. \\
\wp-2,& \text{ if } \vartheta_1=\frac{\wp-3}{2}, \vartheta_1 \equiv 0 \mod 2\\
\wp-1,& \text{ if } \vartheta_1=\frac{\wp-1}{2}, \\
\wp-2,& \text{ if } 1<\vartheta_1<\frac{\wp+1}{2},\vartheta_1 \equiv 0 \mod 2\\
\wp-3,& \text{ if } \frac{\wp+2}{2}<\vartheta_1< \wp-3 ,\vartheta_1 \equiv 0 \mod 2\\
\wp-2,& \text{ if } \frac{\wp-1}{2}<\vartheta_1<\wp-2 ,\\
\wp-2,& \text{ if } \vartheta_1=\wp-2.
\end{array}
\end{cases}
\]
\begin{itemize}
\item If $|\vartheta_1-\vartheta_2|=\wp-1$, then
\end{itemize}
$$
|\mathfrak{R}\{1,\wp\}|= \wp-1.
$$
Noting that the minimum and maximum cardinalities of the resolving sets 
are $\frac{3\wp-5}{4}$ and $\wp-1$, respectively, we conclude by relation 
\eqref{rel} and Theorem~\ref{NT} that
$$
\frac{\wp}{\wp-1}\leq{ldim}_{f}(T_\wp<1, 2, \wp-1>)\leq \frac{4\wp}{3\wp-5}.
$$
The proof is complete.
\end{proof}

\begin{Theorem}
Let $G=T_{2^k}<1,2^{k-2}, 2^{k-1}>$, $k\geq 3$,  
be a Toeplitz graph. The bounds of the local fractional metric dimension of 
$G=T_{2^k}<1, 2^{k-2}, 2^{k-1}>$ is
$$
\frac{2^{k}}{2^{k}-2}\leq{ldim}_{f}(T_\wp<1, 2, \wp-1>)\leq \frac{2^{k}}{8}=2^{k-3}.
$$
\end{Theorem}

\begin{proof} 
The Toeplitz graph $G=T_{2^k}<1,2^{k-2}, 2^{k-1}>$ has order $k\geq 3$. 
Then, we have the following resolving sets:
\begin{itemize}
\item   $\text{ If }   $$|i-h|=1$, then 
$$
\mathfrak{R}\{i,h\}=\mathrm{V}(G)\setminus\left\{i, h\right\}.
$$\item   $\text{ If }   i= 2^{k-3}(\frac{n-2}{2})$, $h=2^{k-3}(\frac{n+2}{2})$, $n=2^k$, then 

\begin{adjustwidth}{-\extralength}{0cm}
$$
\mathfrak{R}\{i,h\}=\mathrm{V}(G)\setminus\left\{1,n\right\}
\bigcup\left\{2,n-1\right\}
\bigcup^{k-2}_{x=2}\left\{2^{x}-1,n-(2^{x}-2)\right\}
\bigcup\left\{\frac{n}{2}\right\}
\bigcup\left\{\frac{n}{2}+2^{k-2}+1\right\}.
$$
\end{adjustwidth}
\item  $\text{ If }  i = 2^{k-3}(\frac{n}{2})$, $h=2^{k-3}(\frac{n+4}{2})$, $n=2^k$, then
$$
\mathfrak{R}\{i,h\}=\mathrm{V}(G)\setminus\{1,n\}
\bigcup\{2,n-1\}
\bigcup\{3,n-2\}
\bigcup\left\{\frac{n}{4}\right\}
\bigcup\left\{\frac{n}{4}+2^{k-2}+1\right\}.
$$
\item \text{ If }  $i= 2^{k-3}+1+a$, $h=2^{k-3}+2+a$, ${0}\leq{a}\leq{2^{k-2}-12}$,
$$
\mathfrak{R}\{i,h\}=\mathrm{V}(G)\setminus\{1+a,2^{k-2}+2+a,2^{k}-2^{k-2}+2+a\}.
$$
\item  $\text{ If }  i=2^{k-1}+2^{k-3}+1+a, h=2^{k-1}+2^{k-3}+2+a,{0}\leq{a}\leq{2^{k-2}-23},$
$$
\mathfrak{R}\{i,h\}=\mathrm{V}(G)\setminus\{1+0,2^{k-1}+1+a,2^{k}-2^{k-2}+2+a\}.
$$
\item  $\text{ If }  i=2^{k-2}+a, h=2^{k-1}+a, ~a=1,2,3,\ldots,2^{k-2},$ then
$$
\mathfrak{R}\{i,h\}=\mathrm{V}(G)\setminus
\bigcup^{2^{k-2}+a-k}_{j=1}\{j\} \bigcup^{a-2^{k-1}+2+1}_{l=0}\{l\}
\bigcup\left\{\frac{2^{k-2}+2^{k-1}+2a}{2}\right\}
\bigcup^{2^{k}}_{j=2^{k-1}+a+k}\{j\}.
$$
\item $\text{ If }  i=2^{k-2}+a, h=2^{k-1}+a, a=1,2,3,\ldots,2^{k-2}$, then
$$
\mathfrak{R}\{i,h\}=\mathrm{V}(G)\setminus\bigcup^{2^{k-2}+a-k}_{j=1}\{j\} 
\bigcup\left\{\frac{2^{k-2}+2^{k-1}+2a}{2}\right\}
\bigcup^{2^{k-2}-1}_{j=0} \{2^{k-2}+a+k+j\}.
$$
\item  $\text{ If }  i=a, h=2^{k-3}+a,a=1,2,3,\ldots,2^{k-2}$, then
$$
\mathfrak{R}\{i,h\}=\mathrm{V}(G)
\setminus\bigcup^{2^{k-2}-2}_{j=0}\{a-2^{k-3}-1+j\}.
$$
\item  $\text{ If }  i=2^{k-3}, h=2^{k-3}+2^{k-2}, $ then
$$
\mathfrak{R}\{i,h\}=\mathrm{V}(G)\setminus\{2^{k-2}\}
\bigcup^{2^{k-2}}_{l=1}\{2^{k-1}+l\}.
$$
\item  $\text{ If }  i=2^{k-3}-j, h=2^{k-3}+2^{k-2}-j, 
j=1,2,3, \ldots ,2^{k-3}-1,$ then
$$
\mathfrak{R}\{i,h\}=\mathrm{V}(G)\setminus\{2^{k-2}-j\}
\bigcup^{2^{k-2}}_{l=1}\{2^{k-1}+l-j\}\bigcup^{i}_{l=1}\{2^{k}+1-l\}.
$$
\item  $\text{ If } i=2^{k-3}+2^{k-1} , h=2^{k-3}+2^{k-2}+2^{k-1},$ then
$$
\mathfrak{R}\{i,h\}=\mathrm{V}(G)\setminus\{2^{k-2}\}
\bigcup^{2^{k-2}-2}_{l=0}\{2^{k-2}+1+l\}\bigcup\{2^{k}-2^{k-2}\}.
$$
\item  $\text{ If } i=2^{k-3}+2^{k-1}-j, 
h=2^{k-3}+2^{k-2}+2^{k-1}-j,j=1,2,3, \ldots ,2^{k-3}-1,$
$$
\mathfrak{R}\{i,h\}=\mathrm{V}(G)\setminus\{2^{k-2}-j\}
\bigcup^{2^{k-2}-2}_{l=0}\{2^{k-2}+1+l\}\bigcup\{2^{k}-2^{k-2}-j\}\\
\bigcup^{j}_{l=1}\{2^{k}+1-l\}.
$$
\item  $\text{ If } i=2^{k-3}+1 , h=2^{k-3}+2^{k-2}+1,$ then
$$
\mathfrak{R}\{i,h\}=\mathrm{V}(G)
\setminus\{2^{k-2}-1\}\bigcup^{2^{k-2}}_{l=1}\{2^{k-1}+1+l\}.
$$
\item  $\text{ If } i=2^{k-3}+1+j, 
h=2^{k-3}+2^{k-2}+1+j,j=1,2,3, 
\ldots ,2^{k-3}-1,$
$$
\mathfrak{R}\{i,h\}=\mathrm{V}(G)\setminus\{2^{k-2}-1-j\}
\bigcup^{2^{k-2}}_{l=1}\{2^{k-1}+1+l-j\}\bigcup^{j}_{l=1}\{l\}.
$$
\item  $\text{ If } i=2^{k-3}+2^{k-1}+1,  
h=2^{k-3}+2^{k-2}+2^{k-1}+1,$ then
$$
\mathfrak{R}\{i,h\}=\mathrm{V}(G)\setminus\{2^{k-2}+1\}
\bigcup^{2^{k-2}-2}_{l=0}\{2^{k-2}+2+l\}.
$$
\item $\text{ If } i=2^{k-3}+2^{k-1}+1+j, 
h=2^{k-3}+2^{k-2}+2^{k-1}+1+j,
j=1,2,3, \ldots ,2^{k-3}-1$, then
$$
\mathfrak{R}\{i,h\}=\mathrm{V}(G)\setminus\{2^{k-2}+1-j\}
\bigcup^{2^{k-2}-2}_{l=0}\{2^{k-2}+2+l\}
\bigcup\{2^{k}-2^{k-2}-j\}\bigcup^{j}_{l=1}\{l\}.
$$
\end{itemize}
The cardinalities of the above resolving set are
\begin{itemize}
\item \text{ If } $|i-h|=1$, then
$$
|\mathfrak{R}\{i,h\}|= 2^{k}-2.
$$

\item $\text{ If } i= 2^{k-3}(\frac{\wp-2}{2}), 
h=2^{k-3}(\frac{n+2}{2}),n=2^k,$ then
$$
|\mathfrak{R}\{i,h\}|= 2^{k}-2k.
$$

\item  $\text{ If }  i = 2^{k-3}(\frac{n}{2})$, 
$h=2^{k-3}(\frac{n+4}{2}),n=2^k$, then
$$
|\mathfrak{R}\{i,h\}|= 8.
$$

\item  $\text{ If }  i= 2^{k-3}+1+a$, $h=2^{k-3}+2+a$, 
${0}\leq{a}\leq{2^{k-2}-12}$,
$$
|\mathfrak{R}\{i,h\}|= 2^{k}-3(2^{k-2}-11).
$$

\item  $\text{ If }  i=2^{k-1}+2^{k-3}+1+a$, 
$h=2^{k-1}+2^{k-3}+2+a,{0}\leq{a}\leq{2^{k-2}-23}$,
$$
|\mathfrak{R}\{i,h\}|= 2^{k}-2(2^{k-2}-22)+1.
$$

\item  $\text{ If }  i=2^{k-2}+a, h=2^{k-1}+a$, $a=1,2,3,\ldots,2^{k-2}$, then
$$
|\mathfrak{R}\{i,h\}|= 2k-2^{k-2}-a-3.
$$

\item  $\text{ If }  i=2^{k-2}+a, h=2^{k-1}+a$, $a=1,2,3,\ldots,2^{k-2}$, then
$$
|\mathfrak{R}\{i,h\}|= 2^{k}-2^{k-1}-k+a+1.
$$

\item  $\text{ If }  i=a, h=2^{k-3}+a,a=1,2,3,\ldots,2^{k-2},$ then
$$
|\mathfrak{R}\{i,h\}|=2^{k}-2^{k-2}+1.
$$

\item  $\text{ If }  i=2^{k-3}, h=2^{k-3}+2^{k-2}, $ then
$$
|\mathfrak{R}\{i,h\}|=2^{k}-2^{k-2}-1.
$$

\item  $\text{ If }  i=2^{k-3}-j$, $h=2^{k-3}+2^{k-2}-j$, 
$j=1,2,3, \ldots ,2^{k-3}-1$, then
$$
|\mathfrak{R}\{i,h\}|= 2^{k}-2^{k-2}-1-i.
$$

\item $\text{ If } i=2^{k-3}+2^{k-1}$, 
$h=2^{k-3}+2^{k-2}+2^{k-1}$, then
$$
|\mathfrak{R}\{i,h\}|= 2^{k}-2^{k-2}-1.
$$

\item $\text{ If } i=2^{k-3}+2^{k-1}-j$, 
$h=2^{k-3}+2^{k-2}+2^{k-1}-j$, 
$j=1,2,3, \ldots ,2^{k-3}-1$, then
$$
|\mathfrak{R}\{i,h\}|= 2^{k}-2^{k-2}-j-1.
$$

\item  $\text{ If } i=2^{k-3}+1$, $h=2^{k-3}+2^{k-2}+1$, then
$$
|\mathfrak{R}\{i,h\}|= 2^{k}-2^{k-2}-1.
$$

\item $\text{ If } i=2^{k-3}+1+j$, 
$h=2^{k-3}+2^{k-2}+1+j$, $j=1,2,3, \ldots ,2^{k-3}-1$, then
$$
|\mathfrak{R}\{i,h\}|= 2^{k}-2^{k-2}-j.
$$

\item  $\text{ If } i=2^{k-3}+2^{k-1}+1$,  
$h=2^{k-3}+2^{k-2}+2^{k-1}+1$, then
$$
|\mathfrak{R}\{i,h\}|= 2^{k}-2^{k-2}.
$$

\item $\text{ If } i=2^{k-3}+2^{k-1}+1+j$, 
$h=2^{k-3}+2^{k-2}+2^{k-1}+1+j$, 
$j=1,2,3, \ldots ,2^{k-3}-1$, then
$$
|\mathfrak{R}\{i,h\}|= 2^{k}-2^{k-2}-j-1.
$$
\end{itemize}
Since the minimum and maximum cardinalities 
of the resolving sets are $8$ and $2^{k}-2$, respectively, 
thus, by relation \eqref{rel} and Theorem~\ref{NT}, we have
$$
\frac{2^{k}}{2^{k}-2}\leq{ldim}_{f}(T_\wp<1, 2, \wp-1>)
\leq \frac{2^{k}}{8}=2^{k-3},
$$
which completes the proof.
\end{proof}

\begin{Theorem}
Let $T_{2p} < 2,p >$ be a Toeplitz graph. 
The bounds of the local fractional metric 
dimension of $T_{2p}< 2,p >$ is
$$
\frac{2p}{2p-1}\leq{ldim}_{f}(T_{2p}<2, p>)
\leq\frac{2p}{2p-2}=\frac{p}{p-1}.
$$
\end{Theorem}

\begin{proof}
Let $T_{2p}< 2,p >$  be a Toeplitz graph. 
One has the following resolving sets:
\[
\mathrm{ }\mathfrak{R}\{\alpha,\beta\}\ 
= 
\begin{cases}
\begin{array}{ll}
T_{2p}< 2,p >(V)/\{p+i-1\},& \text{ if } \alpha=p-2+i, \beta=p+i, i=0,1,2,3,4,\\
T_{2p}< 2,p >(V)/\{p+2\},& \text{ if } \alpha=i, \beta=p+i, i=1,\\
T_{2p}< 2,p >(V)/\{p-1\}.& \text{ if } \alpha=p, \beta=2p,
\end{array}
\end{cases}
\]

\begin{adjustwidth}{-\extralength}{0cm}
\[
\mathfrak{R}\{\alpha,\beta\}\ 
= 
\begin{cases}
\begin{array}{ll}
T_{2p}< 2,p >(V)/\{i-1,p+3\},
& \text{ if } \alpha=i, \beta=p+2, i=2,\\
T_{2p}< 2,p >(V)/\{i+1,p+3+i\},
& \text{ if } \alpha=i, \beta=2+i, i=1,\ldots, p-3,\\
T_{2p}< 2,p >(V)/\{i+1,p+3+i\},
& \text{ if } \alpha=i+2, \beta=p+i+2, i=1,\ldots, p-3,\\
T_{2p}< 2,p >(V)/\{i,p+2+i\},
& \text{ if } \alpha=i+6, \beta=p+i+3, i=1,\ldots, p-3.
\end{array}
\end{cases}
\]
\end{adjustwidth}
The cardinalates of the above resolving sets are
\[
|\mathfrak{R}\{\alpha,\beta\}|\ 
= 
\begin{cases}
\begin{array}{ll}
2p-1,& \text{ if } \alpha=p-2+i, \beta=p+i, i=0,1,2,3,4,\\
2p-1,& \text{ if } \alpha=i, \beta=p+i, i=1,\\
2p-1.& \text{ if } \alpha=p, \beta=2p,
\end{array}
\end{cases}
\]
\[
|\mathfrak{R}\{\alpha,\beta\}|\ 
= 
\begin{cases}
\begin{array}{ll}
2p-2,& \text{ if } \alpha=i, \beta=p+2, i=2,\\
2p-2,& \text{ if } \alpha=i, \beta=2+i, i=1,2,3,\ldots, p-3,\\
2p-2,& \text{ if } \alpha=i+2, \beta=p+i+2, i=1,2,3,\ldots, p-3,\\
2p-2,& \text{ if } \alpha=i+6, \beta=p+i+3, i=1,2,3,\ldots, p-3.
\end{array}
\end{cases}
\]
In general, $2p-2$ is the minimal cardinality and $2p-1$ 
is the maximal cardinality for all prime numbers. So,
by relation (\ref{rel}) and Theorem~\ref{NT}, we have 
$$
\frac{2p}{2p-1}\leq{ldim}_{f}(T_{2p}<2, p>)
\leq\frac{2p}{2p-2}=\frac{p}{p-1}.
$$
The proof is complete.
\end{proof}

\begin{Theorem}
Let $T_{3p}< 3,p>$ be a Toeplitz graph where $p$ is a prime number. 
The local fractional metric dimension of
$T_{3p}< 3,p >$ is equal to one, that is,
$$
{ldim}_{f}\mathrm(T_{3p} < 3,p >)=1.
$$
\end{Theorem}

\begin{proof}
The resolving sets for each adjacent pair of vertices are
$$ 
\mathrm{ }\mathfrak{R}\{\alpha, \beta\}\ 
= \mathrm{V}(T_{3p} < 3,p >) 
\setminus \emptyset = \mathrm{V}(T_{3p} < 3,p >).
$$
So, $3p$ is the minimal and maximal cardinality for all prime numbers and,
by (\ref{rel}) and Theorem \ref{NT}, 
$$
{ldim}_{f}(T_{3p}<3, p>)=1,
$$
which completes the proof.
\end{proof}


\section{LFMD for zero-divisor graphs Over the Set of Zero Divisors of $\mathds{Z}_n$}
\label{sec:04}

In this section, we compute the LFMD for zero-divisor 
graphs over the set $\mathds{Z}_{n}$. These graphs are 
defined as follows: the vertex set is the set of zero 
divisors of $\mathds{Z}_n$ and two vertices are adjacent 
to each other if their product is a zero divisor. 
They are denoted by $G(Z_{n})$.

\begin{Theorem}
\label{mt:02}
Let $G(Z_{n})$  be a zero-divisor graph. An upper bound 
of the local fractional metric dimension of $G(Z_{n})$ is
\[
{ldim}_{f}(G(Z_{n})) 
= 
\begin{cases}
\begin{array}{ll}
1,& \text{ if } n=2p,\\
1,& \text{ if } n=kp,~p>2, k>2,\\
\frac{p ^{k-1}-1}{2},& \text{ if } n=p^{k}, p>3, 
\end{array}
\end{cases}
\]
where $p$ is a prime number.
\end{Theorem}

\begin{proof}
Let $G(Z_{n})$ be a zero-divisor graph. For $n=2p$ the graph $G(Z_{2p})$ 
is bipartite with partition sets 
$$
X=\{p\},~~~~Y=\{\alpha~: ~\gcd(\alpha,2p)\neq1\}\setminus\{p\}.
$$
For illustration, let us take $n=22$. Then, the zero-divisor graph is shown in Figure~\ref{z22}.
\begin{figure}[H]
\includegraphics[scale=0.7]{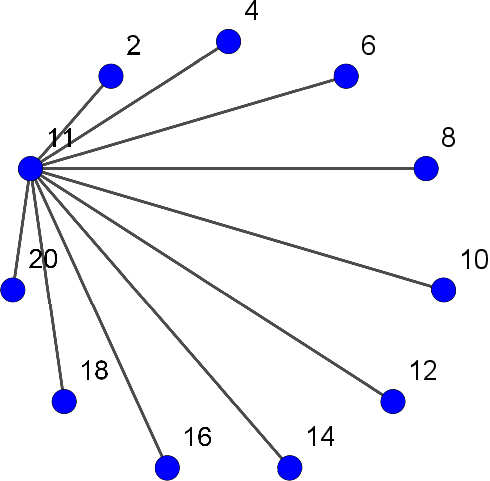}
\caption{The zero-divisor graph $G(Z_{22})$.} \label{z22}
\end{figure}
\noindent By Lemma~\ref{lemma}, ${ldim}_{\vartheta_2}(G(Z_{2p})) = 1$. 

Similarly, for $n=kp$, the graph $G(Z_{kp})$ is also bipartite with partition sets
$$
X=\{p, 2p, \ldots, (k-1)p\},~~~~Y=\{\alpha~: ~\gcd(\alpha,3p)
\neq1\}\setminus\{p, 2p, \ldots, (k-1)p\}.
$$
For observation, we can take $n=33$: 
the zero-divisor graph is shown in Figure~\ref{z33}.
\begin{figure}[H]
\includegraphics[scale=0.7]{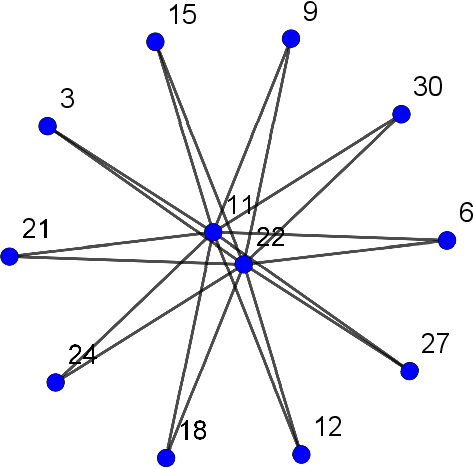}
\caption{The zero-divisor graph $G(Z_{33})$.} \label{z33}
\end{figure}
\noindent By Lemma~\ref{lemma}, ${ldim}_{\vartheta_2}(G(Z_{kp})) = 1$. 

Let $n=p^{k}$ and consider the graph $G(Z_{p^{k}})$. 
The resolving sets for each pair of adjacent vertices are
\begin{equation}
\label{eq:prp:Th6}
\mathfrak{R}\{\alpha,\beta\}=\mathrm{V}({G(Z_{p^{k}})})
\setminus\Big{\{}\mathrm{V}({G(Z_{p^{k}})})\setminus\{\alpha,\beta\}\Big{\}}.
\end{equation}
The cardinality of the set \eqref{eq:prp:Th6} 
is 2 and since the total number of zero divisors 
of $Z_{p^{k}}$ is $p ^{k-1}-1$, we conclude
by Lemma~\ref{lemma} and Theorem~\ref{NT} that 
$$
{ldim}_{f}(G(Z_{kp}))= \frac{p ^{k-1}-1}{2}.
$$
The proof is complete.
\end{proof}


\section{LFMD of zero-divisor graphs Over $\mathds{Z}_{n}\setminus\{0\}$}
\label{sec:05}

\begin{figure}[H]
\includegraphics[scale=0.7]{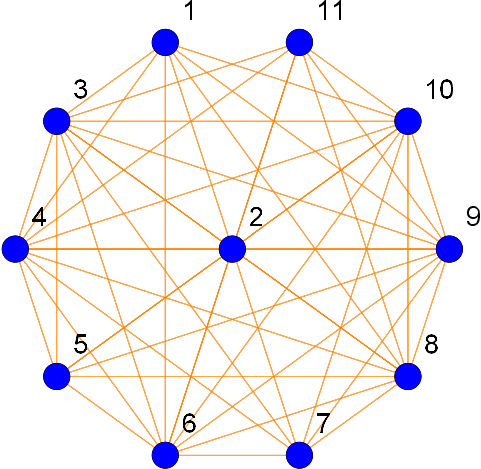}
\caption{The zero-divisor graph $G(\mathds{Z^{*}}_{12})$.} \label{Z12}
\end{figure}

Now, we compute the LFMD for zero-divisor graphs 
over the set $\mathds{Z^{*}}_{n}=\mathds{Z}_{n}\setminus\{0\}$. 
These graphs are defined as follows: the vertex set is 
$\mathds{Z}_{n}\setminus\{0\}=\{1, 2, 3, \ldots, n-1\}$ 
and two vertices are adjacent to each other if their product 
is a zero divisor. They are denoted by $G(\mathds{Z^{*}}_{n})$.
The zero-divisor graph $\mathds{Z}_{12}\setminus\{0\}$ 
is shown in Figure~\ref{Z12}.
The resolving sets corresponding to each pair of adjacent vertices are
\vspace{-12pt}
\begin{adjustwidth}{-\extralength}{0cm}
\begin{align*}
&\mathfrak{R}\{1,2\}=V\{G(\mathds{Z^{*}}_{12})\}\setminus\{3,4,6,8,9,10\},
~~&&\mathfrak{R}\{1,3\}=V\{G(\mathds{Z^{*}}_{12})\}\setminus\{2,4,6,8,9,10\},\\
&\mathfrak{R}\{2,3\}=V\{G(\mathds{Z^{*}}_{12})\}\setminus\{1,4,5,6,7,8,9,10,11\},
~~&&\mathfrak{R}\{1,4\}=V\{G(\mathds{Z^{*}}_{12})\}\setminus\{2,3,6,8,9,10\},\\
&\mathfrak{R}\{4,2\}=V\{G(\mathds{Z^{*}}_{12})\}\setminus\{1,3,5,6,7,8,9,10,11\},
~~&&\mathfrak{R}\{4,3\}=V\{G(\mathds{Z^{*}}_{12})\}\setminus\{1,2,5,6,7,8,9,10,11\},\\
&\mathfrak{R}\{5,2\}=V\{G(\mathds{Z^{*}}_{12})\}\setminus\{3,4,6,8,9,10\}, ~~&&\mathfrak{R}\{5,3\}=V\{G(\mathds{Z^{*}}_{12})\}\setminus\{2,4,6,8,9,10\},\\
&\mathfrak{R}\{5,4\}=V\{G(\mathds{Z^{*}}_{12})\}\setminus\{2,3,6,8,9,10\},
~~&&\mathfrak{R}\{6,1\}=V\{G(\mathds{Z^{*}}_{12})\}\setminus\{2,3,4,8,9,10\},\\
&\mathfrak{R}\{6,2\}=V\{G(\mathds{Z^{*}}_{12})\}\setminus\{1,3,4,5,7,8,9,10,11\}, ~~&&\mathfrak{R}\{6,3\}=V\{G(\mathds{Z^{*}}_{12})\}\setminus\{1,2,4,5,7,8,9,10,11\},\\
&\mathfrak{R}\{6,4\}=V\{G(\mathds{Z^{*}}_{12})\}\setminus\{1,2,3,5,7,8,9,10,11\},
~~&&\mathfrak{R}\{11,9\}=V\{G(\mathds{Z^{*}}_{12})\}\setminus\{2,3,4,6,8,10\},\\
&\mathfrak{R}\{6,5\}=V\{G(\mathds{Z^{*}}_{12})\}\setminus\{2,3,4,8,9,10\},
~~&&\mathfrak{R}\{7,2\}=V\{G(\mathds{Z^{*}}_{12})\}\setminus\{3,4,6,8,9,10\},\\
&\mathfrak{R}\{7,3\}=V\{G(\mathds{Z^{*}}_{12})\}\setminus\{2,4,6,8,9,10\},
~~&&\mathfrak{R}\{7,4\}=V\{G(\mathds{Z^{*}}_{12})\}\setminus\{2,3,6,8,9,10\},\\
&\mathfrak{R}\{7,6\}=V\{G(\mathds{Z^{*}}_{12})\}\setminus\{2,3,4,8,9,10\},
~~&&\mathfrak{R}\{8,1\}=V\{G(\mathds{Z^{*}}_{12})\}\setminus\{2,3,4,6,9,10\},\\
&\mathfrak{R}\{8,2\}=V\{G(\mathds{Z^{*}}_{12})\}\setminus\{1,3,4,5,6,7,9,10,11\},
~~&&\mathfrak{R}\{8,3\}=V\{G(\mathds{Z^{*}}_{12})\}\setminus\{1,2,4,5,6,7,9,10,11\},\\
&\mathfrak{R}\{8,4\}=V\{G(\mathds{Z^{*}}_{12})\}\setminus\{1,2,3,5,6,7,9,10,11\},
~~&&\mathfrak{R}\{8,5\}=V\{G(\mathds{Z^{*}}_{12})\}\setminus\{2,3,4,6,9,10\},\\
&\mathfrak{R}\{8,6\}=V\{G(\mathds{Z^{*}}_{12})\}\setminus\{1,2,3,4,5,7,9,10,11\},
~~&&\mathfrak{R}\{8,7\}=V\{G(\mathds{Z^{*}}_{12})\}\setminus\{2,3,4,6,9,10\},\\
&\mathfrak{R}\{9,1\}=V\{G(\mathds{Z^{*}}_{12})\}\setminus\{2,3,4,6,8,10\},
~~&&\mathfrak{R}\{9,2\}=V\{G(\mathds{Z^{*}}_{12})\}\setminus\{1,3,4,5,6,7,8,10,11\},\\
&\mathfrak{R}\{9,3\}=V\{G(\mathds{Z^{*}}_{12})\}\setminus\{1,2,4,5,6,7,8,10,11\}, ~~&&\mathfrak{R}\{9,4\}=V\{G(\mathds{Z^{*}}_{12})\}\setminus\{1,2,3,5,6,7,8,10,11\},\\
&\mathfrak{R}\{9,5\}=V\{G(\mathds{Z^{*}}_{12})\}\setminus\{2,3,4,6,8,10\}, ~~&&\mathfrak{R}\{9,6\}=V\{G(\mathds{Z^{*}}_{12})\}\setminus\{1,2,3,4,5,7,8,10,11\},\\
&\mathfrak{R}\{9,7\}=V\{G(\mathds{Z^{*}}_{12})\}\setminus\{2,3,4,6,8,10\},
~~&&\mathfrak{R}\{9,8\}=V\{G(\mathds{Z^{*}}_{12})\}\setminus\{1,2,3,4,5,6,7,10,11\},\\
&\mathfrak{R}\{10,1\}=V\{G(\mathds{Z^{*}}_{12})\}\setminus\{2,3,4,6,8,9\},
~~&&\mathfrak{R}\{10,2\}=V\{G(\mathds{Z^{*}}_{12})\}\setminus\{1,3,4,5,6,7,8,9,11\},\\
&\mathfrak{R}\{10,3\}=V\{G(\mathds{Z^{*}}_{12})\}\setminus\{1,2,4,5,6,7,8,9,11\},
~~&&\mathfrak{R}\{10,4\}=V\{G(\mathds{Z^{*}}_{12})\}\setminus\{1,2,3,5,6,7,8,9,11\},\\
&\mathfrak{R}\{10,5\}=V\{G(\mathds{Z^{*}}_{12})\}\setminus\{2,3,4,6,8,9\},
~~&&\mathfrak{R}\{10,6\}=V\{G(\mathds{Z^{*}}_{12})\}\setminus\{1,2,3,4,5,7,8,9,11\},\\
&\mathfrak{R}\{10,7\}=V\{G(\mathds{Z^{*}}_{12})\}\setminus\{2,3,4,6,8,9\},
~~&&\mathfrak{R}\{10,8\}=V\{G(\mathds{Z^{*}}_{12})\}\setminus\{1,2,3,4,5,6,7,9,11\}, \\ &\mathfrak{R}\{10,9\}=V\{G(\mathds{Z^{*}}_{12})\}\setminus\{1,2,3,4,5,6,7,8,11\},
~~&&\mathfrak{R}\{11,2\}=V\{G(\mathds{Z^{*}}_{12})\}\setminus\{3,4,6,8,9,10\},\\
&\mathfrak{R}\{11,3\}=V\{G(\mathds{Z^{*}}_{12})\}\setminus\{2,4,6,8,9,10\},
~~&&\mathfrak{R}\{11,4\}=V\{G(\mathds{Z^{*}}_{12})\}\setminus\{2,3,6,8,9,10\},\\
&\mathfrak{R}\{11,6\}=V\{G(\mathds{Z^{*}}_{12})\}\setminus\{2,3,4,8,9,10\}, ~~&&\mathfrak{R}\{11,8\}=V\{G(\mathds{Z^{*}}_{12})\}\setminus\{2,3,4,6,9,10\},\\
&\mathfrak{R}\{11,10\}=V\{G(\mathds{Z^{*}}_{12})\}\setminus\{2,3,4,6,8,9\}.
\end{align*}
\end{adjustwidth}
The cardinalates of the above resolving set are
\begin{align*}
&|\mathfrak{R}\{1,2\}=5, &&|\mathfrak{R}\{8,2\}|=2,   
&&|\mathfrak{R}\{10,7\}|=5, &&|\mathfrak{R}\{6,3\}|=2, &&|\mathfrak{R}\{9,6\}|=2,\\
&|\mathfrak{R}\{2,3\}|=2, &&|\mathfrak{R}\{8,4\}|=2,  &&|\mathfrak{R}\{10,9\}|=2,   
&&|\mathfrak{R}\{11,9\}|=5,  &&|\mathfrak{R}\{9,8\}|=2,\\
&|\mathfrak{R}\{4,2\}|=2, &&|\mathfrak{R}\{8,6\}|=2, &&|\mathfrak{R}\{11,3\}|=5,   
&&|\mathfrak{R}\{7,2\}|=5,  &&|\mathfrak{R}\{10,2\}|=2,\\
&|\mathfrak{R}\{5,2\}|=5, &&|\mathfrak{R}\{9,1\}|=5, &&|\mathfrak{R}\{11,6\}|=5,   
&&|\mathfrak{R}\{7,4\}|=5, &&|\mathfrak{R}\{10,4\}|=2,\\
&|\mathfrak{R}\{5,4\}|=5, &&|\mathfrak{R}\{9,3\}|=2, &&|\mathfrak{R}\{11,10\}|=5,   
&&|\mathfrak{R}\{8,1\}|=5, &&|\mathfrak{R}\{10,6\}|=2,\\
&|\mathfrak{R}\{6,2\}|=2, &&|\mathfrak{R}\{9,5\}|=5, &&|\mathfrak{R}\{1,3\}|=5,   
&&|\mathfrak{R}\{8,3\}|=2, &&|\mathfrak{R}\{10,8\}|=2,\\
&|\mathfrak{R}\{6,4\}|=2, &&|\mathfrak{R}\{9,7\}|=5, &&|\mathfrak{R}\{1,4\}|=5,  
&&|\mathfrak{R}\{8,5\}|=5,   &&|\mathfrak{R}\{11,2\}|=5.\\
&|\mathfrak{R}\{6,5\}|=5, &&|\mathfrak{R}\{10,1\}|=5, &&|\mathfrak{R}\{4,3\}|=2,  
&&|\mathfrak{R}\{8,7\}|=5,   &&|\mathfrak{R}\{11,4\}|=5,\\
&|\mathfrak{R}\{7,3\}|=5, &&|\mathfrak{R}\{10,3\}|=5, &&|\mathfrak{R}\{5,3\}|=5,  
&&|\mathfrak{R}\{9,2\}|=2,   &&|\mathfrak{R}\{11,8\}|=5,\\
&|\mathfrak{R}\{7,6\}|=5, &&|\mathfrak{R}\{10,5\}|=5, 
&&|\mathfrak{R}\{6,1\}|=5,  &&|\mathfrak{R}\{9,4\}|=2.
\end{align*}
Since the minimum and maximum cardinalities of the resolving 
sets are $2$ and $5$, respectively, we know, by relation 
\eqref{rel} and Theorem~\ref{NT}, that
$$
\frac{11}{5}\leq\mathrm{ldim}_{\mathrm{f}}(G(\mathds{Z^{*}}_{12}))
\leq \frac{11}{2}.
$$

\begin{Theorem}
\label{thm:07}	
Let $G(\mathds{Z^{*}}_{n})$ be a zero-divisors graph over 
$\mathds{Z}_{n}\setminus\{0\}$. The bounds of the local 
fractional metric dimension of
$G(\mathds{Z^{*}}_{n})$ is given by
\[ 
\begin{cases}
\begin{array}{ll}
\frac{2^{k}-1}{2^{k-1}+1}\leq\mathrm{ldim}_{\mathrm{f}}(G(\mathds{Z^{*}}_{2^{k}}))
\leq \frac{2^{k}-1}{2},& \text{ if } n=2^{k}, k>2,\\
\frac{3^{k}-1}{3^{k-1}+1}\leq\mathrm{ldim}_{\mathrm{f}}(G(\mathds{Z^{*}}_{3^{k}}))
\leq \frac{3^{k}-1}{2},& \text{ if } n=3^{k}, k>1,\\
\frac{p^{2}-1}{p^2-p+1}\leq\mathrm{ldim}_{\mathrm{f}}(G(\mathds{Z^{*}}_{p^{2}}))
\leq \frac{p^{2}-1}{2},& \text{ if } n=p^{2}, 
\end{array}
\end{cases}
\]
where $p$ is a prime number.
\end{Theorem}

\begin{proof}
Let $G(\mathds{Z^{*}}_{2^{k}})$ be a zero-divisors graph over 
$\mathds{Z}_{2^{k}}\setminus\{0\}$. We look to the resolving sets 
corresponding to each adjacent vertices.
\begin{itemize}
\item Let $\zeta$ or $\eta$ be an even number, but not both. 
Without any loss of generality, say $\zeta$ is even. Then,
$$
\mathfrak{R}\{\zeta,\eta\}=\mathrm{V}(G(\mathds{Z^{*}}_{2^{k}})
\setminus \Big{\{} \{ \gamma=2t | 1\leq t 
\leq 2^{k-1}-1\}\setminus\{\zeta\} \Big{\}}.
$$
\item If $\zeta$ and $\eta$ are both even numbers, then  
$$
\mathfrak{R}\{\zeta,\eta\}=\mathrm{V}(G(\mathds{Z^{*}}_{2^{k}})
\setminus \Big{\{} \{ t | 1\leq t \leq 2^{k}-1\}
\setminus\{\zeta, \eta\} \Big{\}}.
$$
\end{itemize}
The cardinalities of the above resolving sets are
\begin{itemize}
\item If $\zeta$ or $\eta$ is an even number but not both, then
$$|\mathfrak{R}\{\zeta,\eta\}|=2^{k-1}+1.$$
\item If $\zeta$ and $\eta$ are both even numbers, then
$$
|\mathfrak{R}\{\zeta,\eta\}|=2.
$$
\end{itemize}
Since the minimum and maximum cardinalities of the resolving sets 
are $2$ and $2^{k-1}+1$, respectively, we conclude by relation 
\eqref{rel} and Theorem~\ref{NT} that
$$
\frac{2^{k}-1}{2^{k-1}+1}\leq\mathrm{ldim}_{\mathrm{f}}(G(\mathds{Z^{*}}_{2^{k}}))
\leq \frac{2^{k}-1}{2}.
$$
The proof is similar for $n=3^k$. Let $G(\mathds{Z^{*}}_{p^{2}}\setminus\{0\})$  
be a zero-divisors graph over $\mathds{Z}_{p^{2}}\setminus\{0\}$. 
The resolving sets corresponding to each adjacent vertices are
\begin{itemize}
\item If $\zeta$ or $\eta$ is a multiple of $p$, but not both 
(without any loss of generality, say $\zeta=kp,~~ 1 \leq k \leq p-1$), then
$$
\mathfrak{R}\{\zeta,\eta\}=\mathrm{V}(G(\mathds{Z^{*}}_{p^{2}})
\setminus \Big{\{} \{ \gamma=tp | 1\leq t \leq p-1\}\setminus\{\zeta\} \Big{\}}.
$$
\item If $\zeta$ and $\eta$ are both even numbers, then
$$
\mathfrak{R}\{\zeta,\eta\}=\mathrm{V}(G(\mathds{Z^{*}}_{p^{2}})
\setminus \Big{\{} \{ t | 1\leq t \leq p^{2}-1\}\setminus\{\zeta, \eta\} \Big{\}}.
$$
\end{itemize}
The cardinalities of the above resolving sets are
\begin{itemize}
\item in the first case,
$$
|\mathfrak{R}\{\zeta,\eta\}|=p^{2}-p+1,
$$
\item while in the second
$$
|\mathfrak{R}\{\zeta,\eta\}|=2.
$$
\end{itemize}
Since the minimum and maximum cardinalities of the resolving 
sets are $2$ and $p^2-p+1$, respectively, it follows by relation 
\eqref{rel} and Theorem~\ref{NT} that
$$
\frac{p^{2}-1}{p^2-p+1}
\leq\mathrm{ldim}_{\mathrm{f}}(G(\mathds{Z^{*}}_{p^{2}}))
\leq \frac{p^{2}-1}{2}.
$$
The proof is complete.
\end{proof}


\section{Discussion of Investigated Sequences for the Local Fractional Metric Dimension}
\label{sec:06}

We have examined the asymptotic behavior of families of Toeplitz and zero-divisor graphs 
for the local fractional metric dimension. Note that all obtained sequences 
are constrained. This is illustrated in Table~\ref{t1}.

\begin{table}[H]
\caption{Asymptotic behavior of examined sequences over specific families 
of Toeplitz and zero-divisor graphs for the local fractional metric dimension.}\label{t1}
	\begin{adjustwidth}{-\extralength}{0cm}
		\newcolumntype{C}{>{\centering\arraybackslash}X}
		\begin{tabularx}{\fulllength}{lll}
			\toprule
\boldmath{$G$} & \boldmath{$\mathrm{ldim}_{\mathrm{f}}(G)$} 
& \textbf{Asymptotic Behaviour}\\ \midrule
$T_{\wp}<1, \wp-2>$ &  $ 
\begin{cases}
\begin{array}{ll}
\mathrm{ldim}_{\mathrm{f}}(T_{<1, \wp-2>})=1,
& \text{ if } \wp\in O^{+},\\
\frac{\wp}{\wp-1}\leq\mathrm{ldim}_{\mathrm{f}}(T_{<1, \wp-2>})
\leq\frac{\wp}{\wp-2},&\text{ if } \wp\in E^{+},
\end{array}
\end{cases}$ & Bounded\\ \\
$T_{\wp}<1, 2, \wp-1>$ & $ 
\begin{cases}
\begin{array}{ll}
\mathrm{ldim}_{\mathrm{f}}(T_\wp<1, 2, \wp-1>)=2,& \text{ if } \wp=4,\\
\frac{\wp}{\wp-2}\leq{ldim}_{f}(T_\wp<1, 2, \wp-1>)\leq 2,
&\text{ if }  \wp\equiv 0 \mod 4, \wp \neq 4,\\
\frac{\wp}{\wp-1}\leq{ldim}_{f}(T_\wp<1, 2, \wp-1>)
\leq \frac{2\wp}{\wp+1},&\text{ if }  \wp\equiv 1 \mod 4,\\
1\leq{ldim}_{f}(T_\wp<1, 2, \wp-1>)\leq 2,
&\text{ if }  \wp\equiv 2 \mod 4,\\
\frac{\wp}{\wp-1}\leq{ldim}_{f}(T_\wp<1, 2, \wp-1>)
\leq \frac{4\wp}{3\wp-5},&\text{ if }  \wp\equiv 3 \mod 4.
\end{array}
\end{cases}$ & Bounded\\ \\
$T_{2^{k}}<1, 2^{k-2}, 2^{k-1}>$ & $\frac{2^{k}}{2^{k}-2}
\leq{ldim}_{f}(T_\wp<1, 2, \wp-1>)\leq \frac{2^{k}}{8}=2^{k-3}$ & Unbounded\\ \\
$T_{2p}<2, p>$ & $\frac{2p}{2p-1}\leq{ldim}_{f}(T_{2p}<2, p>)
\leq\frac{2p}{2p-2}=\frac{p}{p-1}$ & Bounded\\ \\
$T_{3p}<3, p>$ & 1 & Constant\\
$G(\mathds{Z}_{2p})$ & 1 & Constant\\ \\
$G(\mathds{Z}_{kp})$ & 1 & Constant\\ \\
$G(\mathds{Z}_{p^{k}})$ & $ \frac{p ^{k-1}-1}{2}$  & Unbounded\\ \\
$G(\mathds{Z}^{*}_{2^{k}})$ & $\frac{2^{k}-1}{2^{k-1}+1}
\leq\mathrm{ldim}_{\mathrm{f}}(G(\mathds{Z^{*}}_{2^{k}}))
\leq \frac{2^{k}-1}{2}$  & Unbounded\\ \\
$G(\mathds{Z}^{*}_{3^{k}})$ & $\frac{3^{k}-1}{3^{k-1}+1}
\leq\mathrm{ldim}_{\mathrm{f}}(G(\mathds{Z^{*}}_{3^{k}}))
\leq \frac{3^{k}-1}{2}$  & Unbounded\\ \\
$G(\mathds{Z}^{*}_{p^{2}})$ & $\frac{p^{2}-1}{p^2-p+1}
\leq\mathrm{ldim}_{\mathrm{f}}(G(\mathds{Z^{*}}_{p^{2}}))
\leq \frac{p^{2}-1}{2}$  & Unbounded\\
			\bottomrule
		\end{tabularx}
	\end{adjustwidth}
\end{table}
The numerical comparison between the computed bounds of 
$T_{2^{k}}<1, 2^{k-2}, 2^{k-1}>$, $G(\mathds{Z}^{*}_{2^{k}})$, 
$G(\mathds{Z}^{*}_{3^{k}})$ and $T_{\wp}<1, \wp-2>$, $T_{\wp}<1, 2, \wp-1>$, 
$\wp\equiv 0 \mod 4$, $T_{\wp}<1, 2, \wp-1>$, $\wp\equiv 1 \mod 4$, 
$T_{\wp}<1, 2, \wp-1>$, $\wp\equiv 2 \mod 4$, 
and $T_{\wp}<1, 2, \wp-1>$, $\wp\equiv 3 \mod 4$ 
is shown in Tables~\ref{T2}, \ref{TT2}, \ref{T3} and \ref{T4}, respectively.

\begin{table}[H]

\caption{Numerical comparison between upper bound values 
of $T_{\wp}<1, \wp-2>$, $T_{\wp}<1, 2, \wp-1>$, $\wp\equiv 0 \mod 4$, 
$T_{\wp}<1, 2, \wp-1>$, $\wp\equiv 1 \mod 4$, $T_{\wp}<1, 2, \wp-1>$, 
$\wp\equiv 2 \mod 4$, and $T_{\wp}<1, 2, \wp-1>$, $\wp\equiv 3 \mod 4$.}
\label{T2}
\newcolumntype{C}{>{\centering\arraybackslash}X}
\begin{tabularx}{\textwidth}{CCCCCC}
\toprule
\boldmath{$\wp$}  &  \boldmath{$T_{\wp}<1, \wp-2>$} & \boldmath{$T_{\wp\equiv0}$} 
&  \boldmath{$T_{\wp\equiv1}$} &  \boldmath{$T_{\wp\equiv2}$}& \boldmath{$T_{\wp\equiv3}$}\\ \hline
8. & 1.14286 & 2 & 1.94118 & 2 & 1.4 \\
9. & 1.125 & 2 & 1.94737 & 2 & 1.39286 \\
10. & 1.11111 & 2 & 1.95238 & 2 & 1.3871 \\
11. & 1.1 & 2 & 1.95652 & 2 & 1.38235 \\
12. & 1.09091 & 2 & 1.96 & 2 & 1.37838 \\
13. & 1.08333 & 2 & 1.96296 & 2 & 1.375 \\
14. & 1.07692 & 2 & 1.96552 & 2 & 1.37209 \\
15. & 1.07143 & 2 & 1.96774 & 2 & 1.36957 \\
16. & 1.06667 & 2 & 1.9697 & 2 & 1.36735 \\
17. & 1.0625 & 2 & 1.97143 & 2 & 1.36538 \\
18. & 1.05882 & 2 & 1.97297 & 2 & 1.36364 \\
19. & 1.05556 & 2 & 1.97436 & 2 & 1.36207 \\
20. & 1.05263 & 2 & 1.97561 & 2 & 1.36066 \\ \bottomrule
\end{tabularx}
\end{table}
\vspace{-12pt}
\begin{table}[H]
\caption{Numerical comparison between lower bound values 
of $T_{\wp}<1, \wp-2>$, $T_{\wp}<1, 2, \wp-1>$, $\wp\equiv 0 \mod 4$, 
$T_{\wp}<1, 2, \wp-1>$, $\wp\equiv 1 \mod 4$, $T_{\wp}<1, 2, \wp-1>$, 
$\wp\equiv 2 \mod 4$, and $T_{\wp}<1, 2, \wp-1>$, $\wp\equiv 3 \mod 4$.}
\label{TT2}
\newcolumntype{C}{>{\centering\arraybackslash}X}
\begin{tabularx}{\textwidth}{CCCCCC}
\toprule
\boldmath{$\wp$}  &  \boldmath{$T_{\wp}<1, \wp-2>$} & \boldmath{$T_{\wp\equiv0}$} 
&  \boldmath{$T_{\wp\equiv1}$} &  \boldmath{$T_{\wp\equiv2}$}& \boldmath{$T_{\wp\equiv3}$}\\ \midrule
8. & 1.06667 & 1.13333 & 1.03125 & 1 & 1.02941 \\
9. & 1.05882 & 1.11765 & 1.02778 & 1 & 1.02632 \\
10. & 1.05263 & 1.10526 & 1.025 & 1 & 1.02381 \\
11. & 1.04762 & 1.09524 & 1.02273 & 1 & 1.02174 \\
12. & 1.04348 & 1.08696 & 1.02083 & 1 & 1.02 \\
13. & 1.04 & 1.08 & 1.01923 & 1 & 1.01852 \\
14. & 1.03704 & 1.07407 & 1.01786 & 1 & 1.01724 \\
15. & 1.03448 & 1.06897 & 1.01667 & 1 & 1.01613 \\
16. & 1.03226 & 1.06452 & 1.01563 & 1 & 1.01515 \\
17. & 1.0303 & 1.06061 & 1.01471 & 1 & 1.01429 \\
18. & 1.02857 & 1.05714 & 1.01389 & 1 & 1.01351 \\
19. & 1.02703 & 1.05405 & 1.01316 & 1 & 1.01282 \\
20. & 1.02564 & 1.05128 & 1.0125 & 1 & 1.0122 \\ \bottomrule
\end{tabularx}
\end{table}
\vspace{-12pt}
\begin{table}[H]
\caption{Numerical comparison between upper bounds values of 
$T_{2^{k}}<1, 2^{k-2}, 2^{k-1}>$, $G(\mathds{Z}^{*}_{2^{k}})$,  
and $G(\mathds{Z}^{*}_{3^{k}})$.}
\label{T3}
\newcolumntype{C}{>{\centering\arraybackslash}X}
\begin{tabularx}{\textwidth}{CCCC}
\toprule
\boldmath{$k$}  &  \boldmath{$T_{2^{k}}<1, 2^{k-2}, 2^{k-1}>$} & \boldmath{$G(\mathds{Z}^{*}_{2^{k}})$} 
&  \boldmath{$G(\mathds{Z}^{*}_{3^{k}})$}\\ \midrule
4. & 2. & 7.5 & 40. \\
5. & 4. & 15.5 & 121. \\
6. & 8. & 31.5 & 364. \\
7. & 16. & 63.5 & 1093. \\
8. & 32. & 127.5 & 3280. \\
9. & 64. & 255.5 & 9841. \\
10. & 128. & 511.5 & 29,524. \\
11. & 256. & 1023.5 & 88,573. \\
12. & 512. & 2047.5 & 265,720. \\
13. & 1024. & 4095.5 & 797,161. \\ \bottomrule

\end{tabularx}
\end{table}

\vspace{-12pt}

\begin{table}[H]\ContinuedFloat
\caption{\textit{Cont}.}
\newcolumntype{C}{>{\centering\arraybackslash}X}
\begin{tabularx}{\textwidth}{CCCC}
\toprule
\boldmath{$k$}  &  \boldmath{$T_{2^{k}}<1, 2^{k-2}, 2^{k-1}>$} & \boldmath{$G(\mathds{Z}^{*}_{2^{k}})$} 
&  \boldmath{$G(\mathds{Z}^{*}_{3^{k}})$}\\ \midrule

14. & 2048. & 8191.5 & $2.39148\times 10^6$ \\
15. & 4096. & 16,383.5 & $7.17445\times 10^6$ \\
16. & 8192. & 32,767.5 & $2.15234\times 10^7$ \\
17. & 16,384. & 65,535.5 & $6.45701\times 10^7$ \\
18. & 32,768. & 131,072. & $1.9371\times 10^8$ \\
19. & 65,536. & 262,144. & $5.81131\times 10^8$ \\
20. & 131,072. & 524,288. & $1.74339\times 10^9$ \\ \bottomrule

\end{tabularx}
\end{table}

\begin{table}[H]
\caption{{Numerical} comparison between lower bounds values of 
$T_{2^{k}}<1, 2^{k-2}, 2^{k-1}>$, $G(\mathds{Z}^{*}_{2^{k}})$,  
and $G(\mathds{Z}^{*}_{3^{k}})$.}
\label{T4}
\newcolumntype{C}{>{\centering\arraybackslash}X}
\begin{tabularx}{\textwidth}{CCCC}
\toprule
\boldmath{$k$}  &  \boldmath{$T_{2^{k}}<1, 2^{k-2}, 2^{k-1}>$} & \boldmath{$G(\mathds{Z}^{*}_{2^{k}})$} 
&  \boldmath{$G(\mathds{Z}^{*}_{3^{k}})$}\\ \midrule
4. & 1.14286 & 1.66667 & 2.85714 \\
5. & 1.06667 & 1.82353 & 2.95122 \\
6. & 1.03226 & 1.90909 & 2.98361 \\
7. & 1.01587 & 1.95385 & 2.99452 \\
8. & 1.00787 & 1.97674 & 2.99817 \\
9. & 1.00392 & 1.98833 & 2.99939 \\
10. & 1.00196 & 1.99415 & 2.9998 \\
11. & 1.00098 & 1.99707 & 2.99993 \\
12. & 1.00049 & 1.99854 & 2.99998 \\
13. & 1.00024 & 1.99927 & 2.99999 \\
14. & 1.00012 & 1.99963 & 3. \\
15. & 1.00006 & 1.99982 & 3. \\
16. & 1.00003 & 1.99991 & 3. \\
17. & 1.00002 & 1.99995 & 3. \\
18. & 1.00001 & 1.99998 & 3. \\
19. & 1. & 1.99999 & 3. \\
20. & 1. & 1.99999 & 3. \\ \bottomrule
\end{tabularx}
\end{table}


\section{Conclusions}
\label{sec:07}

The problem of determining the metric dimension of a given graph for a given integer
is NP-hard. For big graphs, this means that the computation needed to find the metric 
dimension rises exponentially as the graph's size does, making it a difficult task. 
Metric dimension has proven invaluable in addressing problems across various disciplines, 
such as computer science, minimum distance calculations, and chemical graph theory. 
It has played a significant role in improving the representation of networks 
and graphs in real-world applications. Here, we investigated and discussed the asymptotic 
behavior of sequences for the local fractional metric dimension in families of Toeplitz 
and zero-divisor graphs. This sort of study is required in order to understand 
how the attributes of a particular graph or its properties behave as the graph's size increases. 
Without having to manually compute the parameters, asymptotic behavior enables 
researchers to make predictions about how the parameters will behave in bigger networks. 
It is also worth observing that constant, bounded, and unbounded metric 
families all include the obtained upper boundaries for the local fractional 
metric dimension. This implies that the metric dimension and associated 
ideas may be categorized and comprehended within certain families or categories of graphs, 
which can be useful for creating effective algorithms and resolving real-world issues 
in a variety of fields. We trust that our work has several potential applications. 
For instance, cryptography depends on coding and decoding, whereas the most attractive coding 
depends on linear programming problems. If the domain of the coding is the set 
of consecutive integers or some particular class of integers, then its programming 
is very suitable for cryptography. In our proposed work, the class of Toeplitz graphs 
is the zero-divisor graphs that depend on a set of consecutive integers 
and zero divisors of the commutative ring. We claim that our results can be used 
for coding and decoding in cryptography.

For future developments, we mention here some interesting open problems:
\begin{itemize}
\item Toeplitz graph generalization is still up for debate. For such possible
generalizations, it would be great to determine some upper bounds of the local 
fractional metric dimension for $T_nS$ for each subset $S$ 
of $\{1, 2, \ldots, n\}$.
\item For any positive integer $n$, upper bounds for the LFMD of zero-divisor graphs $G(\mathds{Z}_n)$ are yet unknown.
\item For any positive integer $n$, to obtain upper bounds for the LFMD 
of zero-divisor graphs over the set $G(\mathds{Z}^{*}_n)$ 
is still open.
\end{itemize}


\vspace{6pt} 


\authorcontributions{Conceptualization, A.S.A. and S.A.; 
methodology, A.S.A.;
software, M.A.;
validation, A.S.A, S.A. and D.F.M.T.; 
formal analysis, S.A.;
investigation, S.A. and D.F.M.T.;
resources, M.A.;
writing---original draft preparation, S.A. and D.F.M.T.;
writing---review and editing, S.A., M.A. and D.F.M.T.;
visualization, A.S.A.;
supervision, S.A. and D.F.M.T.;
project administration, A.S.A. and D.F.M.T.;
funding acquisition, A.S.A. and D.F.M.T.
All authors have read and agreed to the published version of the manuscript.}


\funding{A.S.A. was supported by Princess Nourah bint Abdulrahman University 
under Researchers Supporting Project PNURSP2023R231, Riyadh, Saudi Arabia;
D.F.M.T. was supported by the Funda\c{c}\~{a}o 
para a Ci\^{e}ncia e a Tecnologia, I.P. (FCT, Funder ID = 50110000187) 
under CIDMA Grants UIDB/04106/2020 and UIDP/04106/2020, Portugal.}


\informedconsent{Not applicable.}

\dataavailability{No data were used to support this study.}

\acknowledgments{A.S.A. is grateful to the support 
of Princess Nourah bint Abdulrahman University 
under Project PNURSP2023R231, Riyadh, Saudi Arabia;
D.F.M.T. to the support by FCT, projects UIDB/04106/2020 and UIDP/04106/2020, Portugal.
The authors would like to express their gratitude 
to four anonymous reviewers for their meticulous review 
of the submitted work, as well as for their valuable questions, 
comments, and suggestions. In light of these insights, the authors 
have made revisions to the manuscript, resulting 
in a noticeable enhancement in its quality.}


\conflictsofinterest{The authors declare no conflict of interest.}


\begin{adjustwidth}{-\extralength}{0cm}
	
\reftitle{References}


\PublishersNote{}

\end{adjustwidth}


\end{document}